\documentclass[letter,12pt]{article}%
\usepackage{amsmath,bm}
\usepackage{amsfonts}
\usepackage{amssymb}
\usepackage{graphicx}
\usepackage{rotating}
\usepackage[height=9in,left=1in,right=0.75in,bottom=1in]{geometry}
\usepackage[breaklinks=true,bookmarksopen=true,colorlinks=true,citecolor=blue]%
{hyperref}
\usepackage[authoryear]{natbib}
\usepackage{color}
\usepackage{colortbl}
\usepackage{paralist}
\usepackage{amsmath}%
\providecommand{\U}[1]{\protect\rule{.1in}{.1in}}
\RequirePackage{hypernat}
\newtheorem{theorem}{Theorem}[section]

\newtheorem{corollary}{Corollary}[section]

\newtheorem{lemma}[theorem]{Lemma}

\newtheorem{proposition}{Proposition}[section]
\newtheorem{remark}{Remark}[section]

\catcode`~=11
\newcommand{\urltilde}{\kern -.15em\lower .7ex\hbox{~}\kern .04em}
\catcode`~=13
\def \@seccntformat#1{\csname the#1\endcsname.\quad}
\makeatother
\hypersetup{pdftitle={Escanciano}, pdfsubject={Semiparametric Identification and Fisher Information}, pdfauthor={Juan Carlos Escanciano}, pdfkeywords={Identification; Irregular Identification; Semiparametric Models; Fisher Information} }
\begin{document}

\title{Semiparametric Identification and Fisher Information\thanks{First version:
September 20th, 2016.}}
\author{Juan Carlos Escanciano\thanks{Department of Economics, Universidad Carlos III
de Madrid, email: jescanci@eco.uc3m.es. Research funded by the Spanish
Programa de Generaci\'{o}n de Conocimiento, reference number
PGC2018-096732-B-I00. I thank Michael Jansson, Ulrich M\"{u}ller, Whitney
Newey, Jack Porter, Pedro Sant'Anna, Ruli Xiao, and seminar participants at
BC, Indiana, MIT, Texas A\&M, UBC, Vanderbilt and participants of the 2018
Conference on Identification in Econometrics for useful comments. }\\\textit{Universidad Carlos III de Madrid}}
\date{September 25th, 2019}
\maketitle

\begin{abstract}
This paper provides a systematic approach to semiparametric identification
that is based on statistical information as a measure of its \textquotedblleft
quality\textquotedblright. Identification can be regular or irregular,
depending on whether the Fisher information for the parameter is positive or
zero, respectively. I first characterize these cases in models with densities
linear in a nonparametric parameter. I then introduce a novel
\textquotedblleft generalized Fisher information\textquotedblright. If
positive, it implies (possibly irregular) identification when other conditions
hold. If zero, it implies impossibility results on rates of estimation. Three
examples illustrate the applicability of the general results. First, I find
necessary conditions for semiparametric regular identification in a structural
model for unemployment duration with two spells and nonparametric
heterogeneity. Second, I show irregular identification of the median
willingness to pay in contingent valuation studies. Finally, I study
identification of the discount factor and average measures of risk aversion in
a nonparametric Euler Equation with nonparametric measurement error in consumption.

\vspace{2mm}

\begin{description}
\item[Keywords:] Identification; Semiparametric Models; Fisher Information.

\item[\emph{JEL classification:}] C14; C31; C33; C35\newpage

\end{description}
\end{abstract}

\section{Introduction}

Nonparametric identification is considered the benchmark for reliable
empirical analysis in economics. Unfortunately, many nonparametric economic
models of interest are unidentified with weak assumptions; see, for example,
discrete choice models with nonparametric unobserved heterogeneity. Yet,
certain interesting aspects of such models might be point-identified by the
same weak assumptions, a situation henceforth referred to as
\textit{semiparametric identification}. Although this observation has long
been recognized in economics (see the early discussion in Hurwicz 1950), no
systematic method is currently available for assessing which aspects (i.e.
functionals) of a nonparametric structural parameter are identified and which
are not. Furthermore, even when nonparametric point-identification holds,
there could be many parameters that are only \textquotedblleft irregularly
identified\textquotedblright\ (in the sense of being identified, but having a
zero Fisher information; see e.g. Chamberlain 1986, Heckman 1990, and Khan and
Tamer 2010). This paper aims to establish general conditions for regular and
irregular \textit{semiparametric} identification (or lack thereof) and to
relate these conditions to the concept of statistical information or the
generalizations proposed herein. The results obtained for irregular
identification\ have important practical implications, as any inferences on
such parameters are expected to be unstable in empirical analysis. In
particular, if a parameter is irregularly identified, then no regular
estimator with a parametric rate of convergence exists (see Chamberlain
1986).\footnote{Notable examples of irregularly identified parameters in
econometrics include densities, regression functions, and their derivatives
evaluated at fixed points of continuous variables; regression discontinuity
parameters, see Cattaneo and Escanciano (2017); binary choice coefficients
under Manski's (1975) conditions, see Chamberlain (1986, 2010); sample
selection models, see Heckman (1990), Andrews and Schafgans (1998) and Goh
(2017); mixed proportional models, see Hahn (1994) and Ridder and Woutersen
(2003); average treatment effects, see Khan and Tamer (2010); or interaction
parameters in triangular systems, see Khan and Nekipelov (2018).}

An important observation for relating identification and information is that
identification depends on both linear and nonlinear effects---see Sargan
(1983) and Chen, Chernozhukov, Lee and Newey (2014)---while statistical
information pertains only to linear effects. To establish a useful link
between the two some structure is thus necessary. Moreover, I show that it is
hard to give sufficient conditions for identification in nonlinear models
allowing for plausibly high levels of irregularity (see Section \ref{Irreg}).
These arguments motivate an initial focus on linear models, i.e. models with
densities that are linear in a nonparametric parameter. In these models a more
complete and transparent analysis of semiparametric identification is
possible, permitting both nonparametric unidentification and high degrees of
irregularity. Specifically, I establish necessary and sufficient conditions
for regular and irregular semiparametric identification. Because many
important economic models are not linear but can be written as linear after
reparametrization, or by fixing some parameters, the results obtained for
linear models are widely applicable.

The analysis of linear models already makes explicit that the separation of
irregular identification from no identification is a rather delicate issue.
The classical Fisher information is not useful when identification is not
regular, because it cannot distinguish between irregular identification and no
identification (it is zero in both cases). This paper introduces a new
\textquotedblleft generalized Fisher information\textquotedblright\ that seems
well-suited for irregular cases. If positive, it implies semiparametric
irregular identification when the classical Fisher information is zero and
other conditions hold. If zero, it implies impossibility results on rates of
convergence for estimators, extending Chamberlain's (1986) impossibility
result to slower rates of convergence than parametric.

The identification results are then extended to semiparametric models that are
nonlinear in the parameter of interest but linear in nuisance parameters.
Examples include commonly used linear and nonlinear panel data models and
structural models of unemployment duration (see e.g. Heckman and Singer 1984a,
1984b), among many others. In this setting, it is possible to present simple
sufficient conditions for identification of the main parameter based on the
generalized Fisher information, allowing for moderate irregularity of the main
parameter, arbitrary irregularity for nuisance parameters, and nonparametric unidentification.

As a general rule, the impossibility results derived in this paper on regular
identification and rates of convergence hold for general linear and nonlinear
models. The sufficient conditions for identification require more structure,
though, because nonlinearities can overwhelm linear effects (cf. Chen et al.
2014). Section \ref{Nonlinear} in the Supplemental Appendix presents
sufficient conditions for semiparametric identification in nonlinear models,
noting that the problem becomes particularly challenging when the model is
nonparametrically unidentified.

I illustrate the usefulness of the theory by deriving new identification
results in the structural model of unemployment duration with two spells and
nonparametric heterogeneity recently proposed by Alvarez, Borovickov\'{a} and
Shimer (2016). These authors first establish that their model is unidentified,
and then discuss a prior sign restriction on the parameters that leads to
nonparametric identification. Complementing their results, I characterize the
identified set without prior sign restrictions and obtain necessary conditions
for regular semiparametric identification. As an implication of my results, I
show that the proportion of individuals at risk of severe long term
unemployment is irregularly identified.

A second example shows irregular identification of the median Willingness to
Pay (WTP) in contingent valuation studies. The median WTP is identified under
weak support conditions, as shown below, and it is an important parameter in
this literature. See Carson and Hanemann (2005) for a detailed survey and
Lewbel, McFadden and Linton (2011) and references therein for semiparametric
estimation of moments of WTP. Using the semiparametric identification strategy
of Lewbel (1997), Khan and Tamer (2010) have provided sufficient conditions
for the mean WTP to be irregularly identified. I revisit the identification of
moments of WTP, and then establish irregular identification of the median WTP,
which is a new result. The operator-based approach followed in the present
paper is very different from that taken in Khan and Tamer (2010), see Section
\ref{WTP} for further discussion.

A third example is a nonlinear model for a consumption-based asset pricing
Euler equation with a nonparametric measurement error in consumption. This
example demonstrates how the results of this paper can be applied to
conditional moment models. Nonparametric and semiparametric treatments of
consumption-based asset pricing models, including Newey and Powell (1988),
Chen and Ludvigson (2009), Escanciano and Hoderlein (2010), Lewbel, Linton and
Srisuma (2011), Chen et al. (2014) and Escanciano et al. (2015), do not deal
with measurement error in consumption. Yet, accounting for measurement error
is vital for empirical studies that use household-level data, as shown in,
e.g., Shapiro (1984), Altonji and Siow (1987), Runkle (1991), and Alan,
Attanasio and Browning (2009). I obtain new primitive conditions for regular
identification of the discount factor and measures of risk aversion under more
general specifications of the marginal utility and the measurement error
mechanism than previously considered. In particular, I show that
identification of the discount factor is more robust to assumptions about the
measurement error than is identification of risk aversion measures.

In summary, this paper provides general semiparametric identification and
impossibility results for regular identification and irregular rates of
estimation, and it shows their utility in some example applications. Inference
on irregular parameters is challenging. In particular, establishing rates of
convergence can be a cumbersome task, but fortunately, the literature has
proposed rate-adaptive estimation and inference methods; see, e.g., Andrews
and Shafgans (1998), Khan and Tamer (2010), Chen and Liao (2014), and Chen and
Pouzo (2015). Rate-adaptive methods are recommended for inference about
irregularly identified parameters.

The question of whether irregular identification holds or not, and if so, to
what degree, is still of first-order importance, because with irregular
identification all estimation methods, including rate-adaptive methods, are
expected to be sensitive to the (unknown) data generating process. This paper
demonstrates that the strength of this sensitivity, and thus the quality of
identification, can be measured by the Fisher information or the
generalizations proposed herein.

The rest of the paper is organized as follows. After a literature review,
Section \ref{Setting} sets the statistical framework and introduces three
examples that will be used throughout the paper. Section \ref{Linear}
characterizes identification in linear models and defines the generalized
Fisher information. Section \ref{SemiparametricModels} analyzes semiparametric
models. Section \ref{Examplesdetail} studies three examples in detail: the
nonparametric unemployment duration model of Alvarez et al. (2016), the median
WTP, and the consumption-based Euler equation with measurement error. Section
\ref{Conclusions} concludes. An Appendix contains proofs of the main results,
and a Supplemental Appendix includes further discussion on identification
conditions for linear and nonlinear models.

\section{Literature Review}

The identification problem has a long history in economics; see the seminal
studies by Koopmans (1949), Hurwicz (1950), Koopmans and Reirsol (1950),
Fisher (1966) and Rothenberg (1971). Bekker and Wansbeek (2001) and Dufour and
Liang (2014) provide more recent contributions as well as a survey of existing
results in parametric settings. Chamberlain (1986) shows that a positive
semiparametric Fisher information is necessary for regular estimation in
semiparametric models. Under the explicit assumption of nonparametric
identification, Van der Vaart (1991) shows the equivalence between a positive
semiparametric information and a differentiability condition that is necessary
for regular estimation. He briefly discusses an \textquotedblleft
intuitive\textquotedblright\ local identification condition, but does not
recognize that this condition may be neither sufficient for identification, as
shown in Chen et al. (2014), nor necessary, as shown in Sargan (1983). Bickel
et al. (1998, Chapter 6) and Ishwaran (1999) present impossibility results on
regularity in some exponential and uniform mixture models. Newey (1990)
provide further impossibility theorems.

There are, of course, many papers reporting sufficient conditions for
identification in specific nonparametric models; see the comprehensive reviews
in Matzkin (2007, 2013) and Lewbel (2016). The closest to my paper is Chen et
al. (2014), who provide sufficient conditions for nonparametric local
identification and for regular semiparametric identification for conditional
moment models. These authors recognize the difficulty of studying
semiparametric irregular identification (see Chen et al. 2014, pg. 796), and
do not analyze that case, which is the focus of this paper. The results of the
present study also help in interpreting their nonparametric identification
conditions in terms of statistical information (see Section \ref{Irreg}). Chen
and Santos (2018) investigate local nonparametric regular overidentification.
Khan and Tamer (2010) show irregular identification in two important examples
and investigate rate-adaptive inference. Chen and Liao (2014) and Chen and
Pouzo (2015) provide general inference results for irregular functionals.
Independently of this paper, Bonhomme (2011) studies regular and irregular
identification and estimation of average marginal effects in nonlinear panel
data models with fixed effects. Also related is the identification analysis of
Severini and Tripathi (2006, 2012) in nonparametric IV. These authors
characterize the set of linear continuous functionals that are regularly and
irregularly identified when the nonparametric structural regression is not
necessarily identified. Santos (2011) investigates semiparametric regular
estimation in the IV setting.

The present paper deals broadly with semiparametric identification in
likelihood and conditional moment models. Furthermore, it follows the
tradition of the seminal work by Rothenberg (1971) in linking the
identification problem to the concept of statistical information (and
generalizations proposed herein), albeit in a nonparametric setting.

\section{Setting and Examples}

\label{Setting}

The data is an independent and identically distributed (iid) sample
$Z_{1},...,Z_{n}$ from a distribution $\mathbb{P}$ that belongs to a class of
probability measures $\mathcal{P}=\{\mathbb{P}_{\lambda}:\lambda\in\Lambda\},$
where $\Lambda$ is a subset of a Hilbert space $(\mathbf{H},\langle\cdot
,\cdot\rangle_{\mathbf{H}})$, with inner product $\langle\cdot,\cdot
\rangle_{\mathbf{H}}$ and norm $\left\Vert \cdot\right\Vert _{\mathbf{H}}.$
For example, in parametric models $\Lambda\subset\mathbb{R}^{m}$ and
$\left\Vert \cdot\right\Vert _{\mathbf{H}}=\left\vert \cdot\right\vert $ is
the Euclidean norm. This paper focusses on nonparametric models where
$\Lambda$ is infinite-dimensional, e.g. a subset of a space of probability
densities. The nonparametric parameter that generates the data is denoted by
$\lambda_{0}\in\Lambda,$ i.e. $\mathbb{P}=\mathbb{P}_{\lambda_{0}}.$ The goal
is to find sufficient and necessary conditions for identification of
$\phi(\lambda)$ at $\phi(\lambda_{0}),$ for a functional $\phi(\lambda
):\Lambda\mapsto\mathbb{R}^{p}$, allowing for the full model $\mathcal{P}$ to
be unidentified at $\lambda_{0}$. That is, the equation $\mathbb{P}_{\lambda
}=\mathbb{P}$ may have more than one solution in $\Lambda$. This setting
includes as a special case semiparametric models where $\lambda=(\theta
,\eta)\in\Lambda=\Theta\times H,$ $\Theta\subset\mathbb{R}^{p}$ and $H$ is a
subset of another Hilbert Space $(\mathcal{H},\langle\cdot,\cdot
\rangle_{\mathcal{H}})$. A leading example of functional is the
finite-dimensional parameter, i.e. $\phi(\theta,\eta)=\theta.$ However, the
setting also includes functionals of the nuisance parameter $\phi
(\lambda)=\chi(\eta),$ where $\chi:H\mapsto\mathbb{R}^{p}$, which, despite the
name, may be of interest. For example, $\eta$ can measure unobserved
heterogeneity, and one might be interested in average marginal effects or
policy counterfactuals that involve averaging across a heterogeneous population.

To introduce the definition of identification, let $f_{\lambda}$ be the
density of $\mathbb{P}_{\lambda}$ with respect to (wrt) a $\sigma$-finite
measure $\mu$. Denote by $\mathcal{B}_{\delta}(\lambda_{0})=\{\lambda
\in\Lambda:\left\Vert \lambda-\lambda_{0}\right\Vert _{\mathbf{H}}<\delta\}$ a
ball of radius $\delta$ around $\lambda_{0}.$\bigskip

\noindent\textbf{Definition (Semiparametric Identification)}: $\phi(\lambda)$
is locally identified in $\mathcal{P}$ at $\phi(\lambda_{0})$ if there exists
$\delta>0$ such that for all $\lambda\in\mathcal{B}_{\delta}(\lambda
_{0}),\ f_{\lambda}=f_{\lambda_{0}}$ $\mu$-almost surely ($\mu$-a.s.) implies
$\phi(\lambda)=\phi(\lambda_{0}).$ If this implication holds for all
$\lambda\in\Lambda$, then $\phi(\lambda)$ is (globally) identified at
$\phi(\lambda_{0})$.\bigskip

To simplify the exposition, I simply write \textquotedblleft$\phi(\lambda
_{0})$ is locally identified\textquotedblright\ rather than \textquotedblleft%
$\phi(\lambda)$ is locally identified in $\mathcal{P}$ at $\phi(\lambda_{0}%
)$\textquotedblright, and if \textquotedblleft locally\textquotedblright\ is
dropped then identification is meant to be global. For parametric models, i.e.
$\Lambda\subset\mathbb{R}^{m},$ Fisher (1966) and Rothenberg (1971) show that,
with sufficient smoothness of the model, non-singularity of the Fisher
information matrix is necessary and sufficient for local identification of
$\lambda_{0}$. In nonparametric models a positive information is not necessary
for identification anymore, and this leads to the classification of
identification in regular and irregular (cf. Khan and Tamer 2010).\bigskip

\noindent\textbf{Definition (Regular and Irregular Semiparametric
Identification)}: $\phi(\lambda_{0})$ is (locally) regularly (respectively,
irregularly) identified if it is (locally) identified and its Fisher
Information is positive (respectively, zero).\bigskip

Identification and regularity/irregularity are separate concepts. So, for
example, the negation of regular identification, which is used extensively
throughout the paper, entails two possibilities: irregular identification or
no identification at all.

The usefulness of the results will be illustrated with several examples. In
all these applications, what distinguishes this paper from others in the
literature is the focus on semiparametric identification and its degree
(regular or irregular), rather than on nonparametric identification and
whether it holds or not. \bigskip

\noindent\textbf{Example 1}: \textit{Unemployment Duration with Heterogeneity}%
. Alvarez, Borovickov\'{a} and Shimer (2016) propose a structural model for
transitions in and out of employment that implies a duration of unemployment
given by the first passage time of a Brownian motion with drift, a random
variable with an inverse Gaussian distribution. The parameters of the inverse
Gaussian distribution are allowed to vary in arbitrary ways to account for
unobserved heterogeneity in workers. These authors investigate nonparametric
identification of the distribution of unobserved heterogeneity, which has a
density $\lambda_{0}$ wrt a $\sigma$-finite measure $\pi,$ when two
unemployment spells $Z_{i}=(t_{i1},t_{i2})$ are observed on the set
$\mathcal{T}^{2},$ $\mathcal{T}\subseteq\lbrack0,\infty)$. The reduced form
parameters $(\alpha,\beta)\in\mathbb{R}\times\lbrack0,\infty)$ are functions
of structural parameters. The distribution of $Z_{i}$ is absolutely continuous
with Lebesgue density $f_{\lambda_{0}}(t_{1},t_{2})$ given, up to a
normalizing constant, by%
\begin{equation}
f_{\lambda_{0}}(t_{1},t_{2})=\int_{\mathbb{R}\times\lbrack0,\infty)}%
\frac{\beta^{2}}{t_{1}^{3/2}t_{2}^{3/2}}e^{-\frac{\left(  \alpha t_{1}%
-\beta\right)  ^{2}}{2t_{1}}-\frac{\left(  \alpha t_{2}-\beta\right)  ^{2}%
}{2t_{2}}}\lambda_{0}(\alpha,\beta)d\pi(\alpha,\beta). \label{densityEx1}%
\end{equation}
Alvarez, Borovickov\'{a} and Shimer (2016) show that $\lambda_{0}$ is
nonparametrically identified up to the sign of $\alpha,$ but do not
investigate semiparametric (regular or irregular) identification, which is the
focus of study here. Specifically, I show that the cdf of $\lambda_{0}$ at a
point will be irregularly identified when $\lambda_{0}$ is identified. This
functional is an important parameter. For example, $\phi(\lambda
_{0})=\mathbb{E}\left[  1\left(  \alpha\leq\alpha_{0}\right)  1\left(
\beta\leq\beta_{0}\right)  \right]  ,$ for a fixed $\alpha_{0}<0<\beta_{0},$
quantifies the proportion of individuals at risk of severe long term
unemployment (an individual with parameters $\alpha$ and $\beta,$ $\alpha
\leq\alpha_{0}$ and $\beta\leq\beta_{0},$ has a probability larger or equal
than $1-\exp(2\alpha_{0}\beta_{0})$ of remaining unemployed forever). This
example is studied in detail in Section \ref{ExampleUnemp}. $\blacktriangle
$\bigskip

\noindent\textbf{Example 2}: \textit{Willingness to Pay}. In contingent
valuation studies one observes $Z_{i}=(Y_{i},V_{i},X_{i}),$ where
$Y_{i}=1\left(  W_{i}\leq V_{i}\right)  ,$ i.e. $Y_{i}=1$ if $W_{i}\leq
V_{i},$ and zero otherwise, $V_{i}$ is a continuous random variable chosen by
the researcher, with known distribution $F_{V}$, and $X_{i}$ a $d-$dimensional
vector of covariates. Here, $W_{i}$ is willingness-to-pay of individual $i$
for a new product or resource, which is an unobserved continuous non-negative
random variable. It is assumed that $W_{i}$ and $V_{i}$ are conditional
independent given $X_{i}.$ The density of $Z_{i},$ with respect to a suitable
measure $\mu$, is%
\[
f_{\lambda_{0}}(y,v,x)=\left[  G_{0}\left(  v,x\right)  \right]  ^{y}\left[
1-G_{0}\left(  v,x\right)  \right]  ^{1-y},
\]
where $G_{0}\left(  v,x\right)  =\mathbb{P}\left[  \left.  W_{i}\leq
v\right\vert X_{i}=x\right]  $ and $\lambda_{0}\left(  v,x\right)  =\partial
G_{0}\left(  v,x\right)  /\partial v.$ Here, one parameter of interest is the
median of the distribution of $W,$%
\[
\phi(\lambda_{0})=Median(W).
\]
Lewbel (1997) and Lewbel, McFadden and Linton (2011) investigate nonparametric
and semiparametric estimation of moments $\phi(\lambda_{0})=\mathbb{E}\left[
r(W_{i})\right]  $ and $\phi(\lambda_{0})=\mathbb{E}\left[  r(W_{i},X)\right]
,$ respectively. Khan and Tamer (2010) show that $\mathbb{E}\left[
W_{i}\right]  $ is irregularly identified when the support of $W$ is
unbounded, and discuss rates of convergence for this functional. There is also
an extensive literature for the related binary choice model when $W_{i}$ has
the representation $W_{i}=\theta_{0}^{\prime}X+\varepsilon_{i}.$ In Section
\ref{WTP} of this paper I study identification of the median WTP. This example
is also useful to illustrate the systematic aspect of the proposed method,
i.e. a single approach can be used for different functionals. $\blacktriangle
$\bigskip

\noindent In the following example, I obtain new regular identification
results by an application of the characterization of regular semiparametric
identification to a conditional moment model. \bigskip

\noindent\textbf{Example 3}: \textit{Consumption-based Asset Pricing Models
with Measurement Error}. Consumption-based asset pricing Euler Equations are
important models in economics. When applied to microeconomic data it is vital
to account for measurement error in consumption, as in%
\[
\mathbb{E}\left[  \left.  \theta_{0}\dot{u}_{0}(C_{t+1}^{\ast})R_{t+1}-\dot
{u}_{0}(C_{t}^{\ast})\right\vert \mathcal{F}_{t}\right]  =0,
\]
where $\theta_{0}$ is the discount factor, $\dot{u}_{0}$ is the marginal
utility of consumption $C_{t}^{\ast}$, $R_{t+1}$ is the gross return of an
asset and $\mathcal{F}_{t}$ denotes the $\sigma$-field generated by the
agent's information set at time $t$. The econometrician observes $C_{t},$
which is a noisy measure of $C_{t}^{\ast},$ following the specification%
\[
C_{t}=m(C_{t}^{\ast},\varepsilon_{t}),
\]
where $m$ is unknown and $\varepsilon_{t}$ is the measurement error. The
primitives of the model, $\lambda_{0},$ are $\theta_{0},$ $\dot{u}_{0},$ $m$
and the distribution of $(C_{t+1}^{\ast},R_{t+1},\varepsilon_{t+1},C_{t}%
^{\ast},\varepsilon_{t})$ given $\mathcal{F}_{t}.$ Of particular interest are
$\theta_{0}$ and $\dot{u}_{0}.$ The observed data is $Z_{i}=(C_{t+1,i}%
,C_{t,i},R_{t+1,i},X_{t,i}),$ for a sample of households, and where $X_{t,i}$
is a vector of household characteristics (e.g. family size) in $\mathcal{F}%
_{t,i}.$ The results of this paper are applied to this example to obtain
sufficient conditions for identification of the discount factor, $\phi
(\lambda_{0})=\theta_{0},$ and the Average Arrow-Pratt coefficient of Absolute
Risk Aversion (AARA), $\phi(\lambda_{0})=\mathbb{E}\left[  \left(
-\partial\dot{u}_{0}(C_{t}^{\ast})/\partial C_{t}^{\ast}\right)  /\dot{u}%
_{0}(C_{t}^{\ast})\right]  $. This example is studied in detail in Section
\ref{ExampleCAPM}. $\blacktriangle$

\section{Linear Nonparametric Models}

\label{Linear}

This section first introduces some notation that will be used throughout the
paper. For a generic measure $\nu,$ let $L_{q}(\nu),$ $q\geq1,$ denote the
Banach space of (equivalence classes of) real-valued measurable functions $h$
such that $\left\Vert h\right\Vert _{q,\nu}:=\left(  \int\left\vert
h\right\vert ^{q}d\nu\right)  ^{1/q}<\infty$ (henceforth I drop the sets of
integration in integrals and the qualification $\nu-$almost surely for
simplicity of notation)$.$ So, for example, a function in $L_{q}(\nu)$ is
discontinuous when there is no continuous function in its equivalence class.
Define the Hilbert space $L_{2}$ of $\mathbb{P}-$square integrable measurable
functions with inner product $\langle h,f\rangle=\int hfd\mathbb{P}$ and norm
$\left\Vert h\right\Vert ^{2}=\langle h,h\rangle$ (I drop the dependence on
$q=2$ and $\nu=\mathbb{P}$ in this case). The set $L_{q}^{0}(\nu)$ ($L_{2}%
^{0})$ is the subspace of zero mean functions in $L_{q}(\nu)$ (respectively,
$L_{2}).$ Henceforth, for a generic linear operator $K:\mathcal{G}%
_{1}\mathcal{\rightarrow G}_{2}$, $\mathcal{N}(K):=\{f\in\mathcal{G}%
_{1}:Kf=0\}$ denotes its kernel. Define $B_{0}=\{b\in\mathbf{H}:\lambda
_{0}+b\in\Lambda\}\ $and let $T(\lambda_{0})$ denote the linear span of
elements in $B_{0}.$ The identification results given below will be relative
to the tangent set $T(\lambda_{0})\ $(i.e, relative to $\Lambda).$

The goal is to relate identification with the concept of statistical
information in a nonparametric setting. To that end, let us consider the
\textit{score operator}, see e.g. Begun, Hall, Huang and Wellner (1983),
which, for linear models defined as in Assumption 1 below, is the operator
$S:T(\lambda_{0})\mapsto L_{2}$ given by%
\begin{equation}
Sb\equiv S_{\lambda_{0}}b:=\frac{f_{\lambda_{0}+b}-f_{\lambda_{0}}}%
{f_{\lambda_{0}}}1(f_{\lambda_{0}}>0). \label{2}%
\end{equation}
Note that under linearity of $S,$ this operator has a unique extension from
$B_{0}$ to $T(\lambda_{0}),$ see Debnath and Mikusinski (2005, pg. 26), and
hence $Sb$ is well-defined when $b\in T(\lambda_{0})\setminus B_{0}.$ More
generally, existence of the score operator is necessary for the classical mean
square differentiability assumption, which means that for every path
$\lambda_{t}\in\Lambda$ with $t^{-1}(\lambda_{t}-\lambda_{0})\rightarrow b\in
T(\lambda_{0})\subset\mathbf{H}$ the following holds,
\begin{equation}
\left\Vert \frac{f_{\lambda_{t}}^{1/2}-f_{\lambda_{0}}^{1/2}}{t}-\frac{1}%
{2}Sbf_{\lambda_{0}}^{1/2}\right\Vert _{2,\mu}\rightarrow0\text{ as
}t\downarrow0. \label{3}%
\end{equation}
The definition of $Sb$ in (\ref{3}) is the most commonly used in the
literature and applies equally to linear and nonlinear models. I will use this
latter definition for nonlinear models, but keep the more natural definition
in (\ref{2}) for linear models. Often the path can be taken of the form
$\lambda_{t}=\lambda_{0}+tb$ and $Sb=\partial\log f_{\lambda_{0}+tb}/\partial
t$ is simply the score associated to the parametric submodel $f_{\lambda
_{0}+tb}$ at the \textquotedblleft truth\textquotedblright\ $t=0$, where,
henceforth, derivatives wrt to $t$ are one-sided and evaluated at zero.

This section investigates identification when the density $f_{\lambda}$ and
the functional $\phi$ are linear.\footnote{For the sake of exposition, I refer
to these as linear, although a more mathematically precise name is affine.} To
simplify the exposition, it is assumed that the functional is a scalar$,$ with
the understanding that all the results below have straightforward extensions
to multivariate functionals$.$\bigskip

\noindent\textbf{Assumption 1}: (i) The map $\dot{\phi}:T(\lambda
_{0})\subseteq\mathbf{H}\mapsto\mathbb{R}$ defined by $\dot{\phi}%
(b)=\phi(\lambda_{0}+b)-\phi(\lambda_{0})$ is linear$;$ (ii) the score
operator $S:T(\lambda_{0})\subseteq\mathbf{H}\mapsto L_{2}$ in (\ref{2}) is
well-defined and linear$;$ (iii) for each $b\in\mathcal{N}(S),$ there exists
$c\equiv c(b)\in\mathbb{R}$, $c\neq0,$ such that $\lambda_{0}+cb\in\Lambda;$
(iv) $\dot{\phi}\ $and $S$ are continuous$.$\bigskip

\noindent Assumption 1(i) holds for the leading example of the
finite-dimensional parameter in a semiparametric model, i.e. $\phi(\theta
,\eta)=\theta.$ Section \ref{Nonlinear} in the Supplemental Appendix relaxes
1(i). Assumption 1(ii) holds in Examples 1 and 2. In other models the
linearity assumption holds after a suitable reparametrization. For example,
consider a simple model of a binary outcome $Y^{\ast}$ that is only observed
when $D=1.$ That is, the available data is $Z=(Y,D,X),$ where $Y=Y^{\ast}D$
and $Y^{\ast}$ is independent of $D$ given $X.$ Define $q(x)=\mathbb{E}\left[
\left.  Y^{\ast}\right\vert X\right]  $ and $p(x)=\mathbb{E}\left[  \left.
D\right\vert X\right]  .$ The density of $(Y,D,X)$ is a nonlinear function of
$p$ and $q,$ since for example $\mathbb{P}\left(  \left.  Y=1,D=1\right\vert
X=x\right)  =q(x)p(x),$ but if we reparametrize the density in terms of
$\lambda_{0}=(\lambda_{01},\lambda_{02})$ with $\lambda_{01}(x)=q(x)p(x)$ and
$\lambda_{02}(x)=p(x),$ then the density becomes linear in $\lambda_{0}$.
Section 5 below and Section \ref{Nonlinear} in the Supplemental Appendix relax
Assumption 1(ii). Assumption 1(iii) is only used to prove the necessity of the
main identification condition below, and it can be dropped altogether by
restricting attention only to $b\in B_{0}$, see Remark 2 below, although it
facilitates exposition. Overall, Assumption 1 is convenient because, with this
assumption, identification can be fully characterized. Thus, the
identification results under Assumption 1 provide a benchmark for what can be
achieved in more complicated situations.

Assumption 1(iv) guarantees the existence of the adjoint operator $S^{\ast
}:L_{2}\mapsto T(\lambda_{0})$ of $S,$ satisfying $\langle g,Sb\rangle=\langle
S^{\ast}g,b\rangle_{\mathbf{H}}$ for all $g\in L_{2}$ and $b\in T(\lambda
_{0}).$ The following definition extends the Fisher information matrix to a
nonparametric context (cf. Ko\v{s}evnik and Levit 1976).\bigskip

\noindent\textbf{Definition (Fisher Information)}: The information operator is
defined as $I_{\lambda_{0}}:=S^{\ast}S$.\bigskip

Roughly, $I_{\lambda_{0}}b$ measures the Fisher information of $\lambda_{0}$
in the direction $b\in T(\lambda_{0}),$ i.e. the classical Fisher information
corresponding to $f_{\lambda_{0}+tb}$ at $t=0.$ To establish a link between
$I_{\lambda_{0}}$ and identification, note that, under Assumption 1,
semiparametric identification of $\phi(\lambda_{0})$ will hold if, for all
$b\in T(\lambda_{0}),$%
\[
f_{\lambda_{0}+b}-f_{\lambda_{0}}\equiv f_{\lambda_{0}}1(f_{\lambda_{0}%
}>0)Sb=0\Longrightarrow\phi(\lambda_{0}+b)-\phi(\lambda_{0})\equiv\dot{\phi
}(b)=0,
\]
or, since $f_{\lambda_{0}}>0$ $\mathbb{P}-$a.s. and $\mathcal{N}%
(S)=\mathcal{N}(I_{\lambda_{0}}),$
\begin{equation}
\mathcal{N}(I_{\lambda_{0}})\subset\mathcal{N}(\dot{\phi}). \label{1}%
\end{equation}
The following proposition proves that (\ref{1}), which involves the
nonparametric Fisher information, is necessary and sufficient for
semiparametric identification of $\phi(\lambda_{0})$ under Assumption 1.

\begin{proposition}
\label{Identification}Under Assumption 1, identification of $\phi(\lambda
_{0})$ holds iff (\ref{1}) holds.
\end{proposition}

\noindent Without Assumption 1(ii) both implications of Proposition
\ref{Identification} fail, which motivates the initial focus on linear models.
That (\ref{1}) is not sufficient for identification follows from a
counterexample given in Chen et al. (2014), while that it is not necessary by
Sargan (1983). Some structure is thus needed for the intuitive identifiability
condition (\ref{1}) to be useful for identification. Assumption 1 is a natural
starting point, because with this assumption identification is characterized.
\bigskip

\noindent\textbf{Remark 1}: (i) When $\Lambda$ is a subset of densities wrt a
$\sigma$-finite measure $\pi,$ it is natural to define the score operator as
an operator from $L_{1}(\pi)$ to $L_{1}(\mathbb{P})$ since $\lambda\in
L_{1}(\pi)$ and $Sb\in L_{1}(\mathbb{P})$ without additional assumptions.
Proposition \ref{Identification} holds with these extended definitions as
well. (ii) If a Hilbert space approach is preferred, it is convenient to
introduce directions $b\in L_{2}^{0}(G_{0}),$ where $G_{0}$ is the measure
associated to $\lambda_{0},$ and re-define the score operator in (\ref{2}) and
the functional $\dot{\phi}$ with $b$ replaced by $\lambda_{0}b$ (and
$T(\lambda_{0})$ as the linear span of $B_{0}=\{b\in L_{2}^{0}(G_{0}%
):\lambda_{0}+\lambda_{0}b\in\Lambda\}$). The specific definition used for the
operator $S$ and functional $\dot{\phi}$ will be clear from the context.
\bigskip

\noindent\textbf{Remark 2}: Assumption 1(iii) can be dropped if the
identification condition (\ref{1}) is replaced by%
\begin{equation}
\mathcal{N}(S)\cap B_{0}\subset\mathcal{N}(\dot{\phi})\cap B_{0}. \label{1.1}%
\end{equation}
\bigskip

\noindent The identification condition (\ref{1}) is based on the nonparametric
Fisher information and, as such, it is not useful if the goal is to
disentangle regular and irregular semiparametric identification. To introduce
a more useful characterization, I first define the semiparametric Fisher
information for $\phi$. The information for estimating the parameter
$\psi(t)=\phi(\lambda_{t}),$ $\lambda_{t}:=\lambda_{0}+tb,$ under the density
$f_{\lambda_{t}}$ at $t=0$ is, by the delta method, equal to $\left[
\partial\psi(t)/\partial t\right]  ^{-1}||Sb||^{2}\left[  \partial
\psi(t)/\partial t\right]  ^{-1}=||Sb||^{2}/[\dot{\phi}(b)]^{2}.$ The
semiparametric Fisher information is the infimum of the informations over all
such parametric submodels (cf. Stein 1956) and is given by
\begin{equation}
I_{\phi}=\inf_{b\in\mathcal{B}_{\phi}}\frac{||Sb||^{2}}{\left[  \dot{\phi
}(b)\right]  ^{2}}, \label{Fisher}%
\end{equation}
where $\mathcal{B}_{\phi}:=\left\{  b\in T(\lambda_{0}):\dot{\phi}%
(b)\neq0,\left\vert \dot{\phi}(b)\right\vert \leq1\right\}  .$

By the continuity in Assumption 1(iv), $S$ and $\dot{\phi}$ are uniquely
extended to $\overline{T(\lambda_{0})},$ where henceforth, for a subspace $V,$
$\overline{V}$ denotes the closure of $V$ in the norm topology. Moreover,
there exists an $r_{\phi}\in\overline{T(\lambda_{0})},$ called the Riesz's
representer of $\dot{\phi},$ such that for all $b\in\overline{T(\lambda_{0}%
)},$
\[
\dot{\phi}(b)=\langle b,r_{\phi}\rangle_{\mathbf{H}}.
\]
I can then identify $\dot{\phi}\ $with $r_{\phi},$ and provide identification
results in terms of $r_{\phi}$ using duality. I will make clear in the text
when $S$ is viewed as acting on $T(\lambda_{0})$ or on its extension
$\overline{T(\lambda_{0})}.$ Henceforth, whenever identification is discussed
in terms of $r_{\phi},$ like in the next theorem, I mean $S$ and $\dot{\phi}$
to be defined on $\overline{T(\lambda_{0})}$ (so, for example, $r_{\phi}$ is
guaranteed to exist by the Riesz representation theorem). Let $\mathcal{R}%
(S^{\ast}):=\{f\in\overline{T(\lambda_{0})}:\exists g\in L_{2},S^{\ast}g=f\}$.
Then, the following result provides a full characterization of semiparametric
regular and irregular identification in linear models.

\begin{theorem}
\label{MainTheoremS}Under Assumption 1: (i) $\phi(\lambda_{0})\ $is regularly
identified iff $r_{\phi}\in\mathcal{R}(S^{\ast})$; (ii) $\phi(\lambda_{0}%
)\ $is irregularly identified iff $r_{\phi}\in\overline{\mathcal{R}(S^{\ast}%
)}\setminus\mathcal{R}(S^{\ast});$ and (iii) $\phi(\lambda_{0})\ $is
unidentified iff $r_{\phi}\notin\overline{\mathcal{R}(S^{\ast})}.$
\end{theorem}

\noindent The results that seem to be novel here are the \textquotedblleft
identification\textquotedblright\ parts and the separation of irregular
identification from no identification. The fact that $r_{\phi}\in
\mathcal{R}(S^{\ast})$ is equivalent to $I_{\phi}>0$ is due to Van der Vaart
(1991, Theorem 4.1)$.$ As discussed earlier, the characterizations of Theorem
\ref{MainTheoremS} do not necessarily hold without Assumption 1, and some care
must be exercised to extrapolate these results to nonlinear models. Section
\ref{counterexample} in the Supplemental Appendix presents a nonlinear
counterexample, building on that given in Chen et al. (2014), that shows that
$I_{\phi}>0$ may hold while $\phi(\lambda_{0})$ is not identified.\bigskip

\noindent\textbf{Remark 3}: Assumption 1(iii) can be dropped for part (i) of
Theorem \ref{MainTheoremS} because a positive Fisher information implies
$r_{\phi}\in\mathcal{R}(S^{\ast})$.

\subsection{Irregular Identification and Generalized Information}

\label{Irreg}

When the Fisher information $I_{\phi}$ is zero, it does not provide
information on identification (it cannot distinguish between irregular
identification and no identification). This section introduces a
\textquotedblleft generalized Fisher information\textquotedblright\ that
extends, in a sense described later, the classical Fisher information to
irregular cases, and which is given by%
\begin{equation}
I_{\phi,\rho}=\inf_{b\in\mathcal{B}_{\phi}}\frac{||Sb||^{2}}{\left\vert
\dot{\phi}(b)\right\vert ^{2\rho}}, \label{genFisher}%
\end{equation}
where $1\leq\rho<\infty.$ The classical Fisher information corresponds to
$\rho=1,$ i.e. $I_{\phi}\equiv I_{\phi,1}$. Furthermore, it can be shown that
$I_{\phi,1}\leq I_{\phi,\rho}$ for $1<\rho<\infty$.

The sense in which the generalized Fisher information provides a
generalization of the classical Fisher information is shown in the next two
results. The first result extends the sufficient condition for identification
in Theorem \ref{MainTheoremS}(i) to the irregular case. Under Assumption 1,
$I_{\phi,1}=0$ and $I_{\phi,\rho}>0$ for some $\rho>1$ corresponds to
irregular identification$.$ Because $\left\vert \dot{\phi}(b)\right\vert \leq
I_{\phi,\rho}^{-2\rho}||Sb||^{1/\rho},$ for $b\in\mathcal{B}_{\phi},$ an
interpretation of a positive generalized information is that $\dot{\phi}(b)$
is continuous in $\mathcal{B}_{\phi}$ wrt the Fisher semi-norm $||Sb||,$ with
a modulus of continuity quantified by $1/\rho$ (smaller $\rho$ corresponding
to more regularity). The inequality above directly gives identification on a
restricted set, and the next result proves that this set can be extended to
the whole parameter space $\Lambda$ under Assumption 1.

\begin{theorem}
\label{TheoremR}Let Assumption 1(i-ii) hold. If $I_{\phi,\rho}>0,$ $1\leq
\rho<\infty,$ then $\phi(\lambda_{0})\ $is identified.
\end{theorem}

\noindent As discussed at the end of this section, there is an extensive
literature on nonparametric identification that requires, implicitly or
explicitly, the condition $I_{\phi,\rho}>0$ to establish rates of estimation.
Theorem \ref{TheoremR} shows $I_{\phi,\rho}>0$ is sufficient for
identification, while the next result shows it is necessary for achieving
certain rates of convergence. Let $\mathcal{A}$ denote a class of sequences in
$\Lambda,$ and let $\mathbb{P}_{\lambda}^{n}$ denote the $n-$fold probability
of $\mathbb{P}_{\lambda}$. I provide a formal definition of rate of
convergence (see e.g. Ishwaran 1996, Definition 7).\bigskip

\noindent\textbf{Definition (Rate of Convergence)}: The estimator $T_{n}$ has
a rate of convergence $r_{n}$ on $\mathcal{A}$ for estimating $\phi
(\lambda_{0})$ if for each $\varepsilon>0$ there exists $K(\varepsilon)>0$
such that for each $\{\lambda_{n}\}\in\mathcal{A},$%
\[
\lim\sup_{n\rightarrow\infty}\mathbb{P}_{\lambda_{n}}^{n}\left(  \left\vert
T_{n}-\phi(\lambda_{n})\right\vert >K(\varepsilon)r_{n}\right)  <\varepsilon.
\]

\begin{theorem}
\label{Impossibility}Suppose that for each $\varepsilon>0$ there exists a path
$\lambda_{t}\in\Lambda$ passing through $\lambda_{0}$ such that for all $t$
sufficiently small the following holds: (i) $Ct\leq\left\vert \phi(\lambda
_{t})-\phi(\lambda_{0})\right\vert \leq1$, for some $C>0,$ and (ii)
$\left\Vert \left(  f_{\lambda_{t}}-f_{\lambda_{0}}\right)  /f_{\lambda_{0}%
}\right\Vert ^{2}<\varepsilon t^{2\rho},$ for $1\leq\rho<\infty.$ Suppose
$\mathcal{A}$ contains all sequences $\{\lambda_{n}\}$ for which $\phi
(\lambda_{n})=\phi(\lambda_{0})+O(n^{-1/2\rho}).$ Then, the rate of
convergence for any estimator of $\phi(\lambda_{0})$ on $\mathcal{A}$ must be
slower than $O(n^{-1/2\rho}).$
\end{theorem}

\noindent Conditions (i) and (ii) in Theorem \ref{Impossibility} correspond to
$I_{\phi,\rho}=0,$ since under these conditions $b_{t}\equiv(\lambda
_{t}-\lambda_{0})\in\mathcal{B}_{\phi}$ and%
\[
I_{\phi,\rho}=\inf_{b\in\mathcal{B}_{\phi}}\frac{||Sb||^{2}}{[\dot{\phi
}(b)]^{2\rho}}\leq\frac{\varepsilon}{C^{2\rho}},
\]
and because $\varepsilon>0\ $is arbitrary, it must hold that $I_{\phi,\rho
}=0.$ I stress that Assumption 1 is not required here, so Theorem
\ref{Impossibility} holds for linear and nonlinear models/parameters. Theorem
\ref{Impossibility} extends the impossibility result of Chamberlain (1986) to
the irregular case $\rho>1.$ Its proof relies on techniques, originally due to
Lecam (1973), that bound the total variation distance between the distribution
under $f_{\lambda_{n}}$ and that under $f_{\lambda_{0}},$ for a suitable
sequence $\lambda_{n}.$ Previously, Donoho and Liu (1987) have used similar
techniques based on Hellinger distance to provide lower bounds for convergence
of estimators. The modest contribution of Theorem \ref{Impossibility} is to
connect the zero generalized Fisher information introduced in this paper with
these lower bounds, extending important work by Chamberlain (1986) to the
irregular case.

It is important to stress that there is nothing that prevents the possibility
that $I_{\phi,\rho}=0$ for all $\rho\geq1$, which, under Theorem
\ref{Impossibility}, implies impossibility of polynomial rates. Indeed, in
several important models logarithmic rates are common, see, e.g., Fan (1991)
for classical measurement error problems. This possibility suggests that the
definition of the generalized information and the conditions of Theorem
\ref{Impossibility} should be modified to accommodate severely irregular
cases. For a function $\psi$ that is increasing, non-negative, right
continuous at $0$ and with $\psi(0)\downarrow0,$ one can define the
generalized Fisher information%
\[
I_{\phi,\psi}=\inf_{b\in\mathcal{B}_{\phi}}\frac{||Sb||^{2}}{\psi\left(
\lbrack\dot{\phi}(b)]^{2}\right)  }.
\]
With this modification different degrees of irregularity, including severe
irregularity, are allowed. The case $\psi\left(  \epsilon\right)
=\epsilon^{\rho},$ $1<\rho\leq2,$ corresponds to mild or moderate
irregularity, while $\psi\left(  \epsilon\right)  =\exp(\epsilon)-1$ or
$\psi\left(  \epsilon\right)  =\exp(-1/(\epsilon^{a})),$ with $a>0,$ is
suitable for severe irregularity with possibility of logarithmic rates. A
version of Theorem \ref{Impossibility} that allows for severe irregularity
follows \textit{mutatis mutandis}, simply replacing $t^{2\rho}$ by
$\psi\left(  t^{2}\right)  .$

At this point, it is useful to compare the results of this paper with the
general nonparametric local identification results in Chen et al. (2014) for
conditional moment restrictions models.\ These authors obtain sufficient
conditions for nonparametric identification of linear and nonlinear
conditional moments by suitably restricting the parameter space. When
conditional moments are only Frechet differentiable, they consider the
parameter space to have tangents in $\{b:||\dot{m}b||^{2}>C\left\Vert
b\right\Vert _{\mathbf{H}}^{2}\}$, for the derivative $\dot{m}$ of a
conditional mean operator $m$ and a positive constant $C.$ In our setting, an
analog that allows comparison with statistical information would be
$m(\lambda)=\left(  f_{\lambda}-f_{\lambda_{0}}\right)  /f_{\lambda_{0}},$
with derivative $\dot{m}(\lambda)=S(\lambda-\lambda_{0}).$\footnote{More
precisely, the effective score operator in their conditional moment
restriction model is $S(\lambda-\lambda_{0})=\Sigma^{-1}\dot{m}(\lambda
-\lambda_{0})$ for a conditional variance $\Sigma.$ Under the standard
assumption that the eigenvalues of $\Sigma$ are bounded away from zero, there
is a one to one relation between identification in conditional moment models
as in Chen et al. (2014) and identification in terms of statistical
information as considered here; see e.g. Chamberlain (1992).} On the parameter
space with tangents $\{b:||Sb||^{2}>C\left\Vert b\right\Vert _{\mathbf{H}}%
^{2}\}$ the nonparametric information is positive (i.e. regular nonparametric
identification), which implies a positive semiparametric Fisher information
for \textit{all} linear continuous functionals. Chen et al. (2014) also
consider conditions corresponding to higher order differentiability and these
conditions do allow for irregular semiparametric identification. In their
general case, they restrict tangents to the set $\{b:||Sb||^{2}>C\left\Vert
b\right\Vert _{\mathbf{H}}^{2\rho}\},$ for $\rho>1,$ which implies a positive
generalized Fisher information $I_{\phi,\rho}$ for \textit{all} continuous
linear functionals $\phi$. To prove such result, use $\left\vert \dot{\phi
}(b)\right\vert \leq\left\Vert r_{\phi}\right\Vert _{\mathbf{H}}\left\Vert
b\right\Vert _{\mathbf{H}}\ $and their assumption $\mathcal{B}_{\phi}%
\subset\{b:||Sb||^{2}>C\left\Vert b\right\Vert _{\mathbf{H}}^{2\rho}\}$ to
bound
\[
I_{\phi,\rho}=\inf_{b\in\mathcal{B}_{\phi}}\frac{||Sb||^{2}}{[\dot{\phi
}(b)]^{2\rho}}\geq\frac{C}{\left\Vert r_{\phi}\right\Vert _{\mathbf{H}}%
^{2\rho}}>0.
\]
This shows that the restrictions on neighborhoods in Chen et al. (2014) have a
statistical interpretation in terms of the generalized Fisher information
introduced here (for all continuous linear functionals). The parameter $\rho$
is also linked to the nonlinearity permitted in the model (see Assumption 2 in
Chen et al. 2014 or Assumption 2 below), which typically restricts its values
to $1\leq\rho\leq2.$ The case $\rho>2$ only holds when the second derivative
of $m(\lambda)$ is zero. From this discussion, it seems that it may be hard to
accomodate nonlinear cases for severely ill-posed problems where $I_{\phi
,\rho}=0$ for all $\rho\geq1$ and \textit{some} functional $\phi$. Examples of
such severe irregularity for certain functionals $\phi$ include nonparametric
IV models or measurement error models with Gaussian errors.

To give sufficient conditions for irregular identification that allow for a
variety of degrees of irregularity for different functionals consider adding
the following condition to Assumption 1:\bigskip

\noindent\textbf{Assumption 1}: (v) the score operator $S$ is compact.
\bigskip

\noindent Assumption 1(v) guarantees the existence of a sequence
$\{\lambda_{j},\varphi_{j},\psi_{j}\}_{j=1}^{\infty}$ such that (cf. Kress,
1999, Theorem 15.16) \textit{\ }%
\begin{equation}
S\varphi_{j}=\lambda_{j}\psi_{j}\qquad\text{and}\qquad S^{\ast}\psi
_{j}=\lambda_{j}\varphi_{j}. \label{rel}%
\end{equation}
This is the so called Singular Value Decomposition (SVD) of $S.$ The elements
$\{\varphi_{j}\}_{j=1}^{\infty}$ and $\{\psi_{j}\}_{j=1}^{\infty}$ are
complete orthonormal bases for $\overline{\mathcal{R}(S^{\ast})}$ and
$\overline{\mathcal{R}(S)},$ respectively, and the singular values
$\lambda_{j}$ are the squared-root eigenvalues of the information operator
$I_{\lambda_{0}}=S^{\ast}S:\overline{T(\lambda_{0})}\mapsto\overline
{T(\lambda_{0})}.$ Furthermore, defining for $\beta\in\mathbb{R},$
\[
\mathcal{M}_{\beta}:=\left\{  b\in\overline{T(\lambda_{0})}\text{ such that
}\left\Vert b\right\Vert _{\beta}^{2}:=\sum_{j=1}^{\infty}\lambda_{j}%
^{-2\beta}\left\langle b,\varphi_{j}\right\rangle _{_{\mathbf{H}}}^{2}%
<\infty\right\}  ,
\]
it is well known (see e.g. Carrasco, Florens and Renault 2007) that%
\[
\overline{\mathcal{R}(S^{\ast})}\equiv\mathcal{M}_{0}=\left\{  b\in
\overline{T(\lambda_{0})}\text{ such that }\sum_{j=1}^{\infty}\left\langle
b,\varphi_{j}\right\rangle _{_{\mathbf{H}}}^{2}<\infty\right\}  ,
\]
whereas%
\[
\mathcal{R}(S^{\ast})\equiv\mathcal{M}_{1}=\left\{  b\in\overline
{T(\lambda_{0})}\text{ such that }\sum_{j=1}^{\infty}\lambda_{j}%
^{-2}\left\langle b,\varphi_{j}\right\rangle _{_{\mathbf{H}}}^{2}%
<\infty\right\}  .
\]
The following result gives sufficient, and in some cases necessary, conditions
for regular and irregular identification in terms of $\left\Vert r_{\phi
}\right\Vert _{\beta}.$

\begin{theorem}
\label{ThmSVD}Let Assumption 1 hold. Then (i) $\phi(\lambda_{0})\ $is
regularly identified iff $\left\Vert r_{\phi}\right\Vert _{1}<\infty,$ and in
that case $I_{\phi,1}=\left\Vert r_{\phi}\right\Vert _{1}^{-2}$; (ii) If
$\left\Vert r_{\phi}\right\Vert _{1}=\infty$ but $\left\Vert r_{\phi
}\right\Vert _{\beta}<\infty$ for some $0<\beta<1,$ then%
\[
\inf_{b\in\mathcal{B}_{\phi},||b||_{\mathbf{H}}\leq1}\frac{||Sb||^{2}}%
{[\dot{\phi}(b)]^{2\rho}}\geq\frac{1}{\left\Vert r_{\phi}\right\Vert _{\beta
}^{2\rho}}>0,
\]
for $\rho=1/\beta,$ so $\phi(\lambda_{0})\ $is locally irregularly identified.
\end{theorem}

In part (i) it is possible to give another expression for the efficiency bound
as a variance. The condition $\left\Vert r_{\phi}\right\Vert _{1}<\infty$
implies there exists $i_{\phi}\in L_{2}$ such that $S^{\ast}i_{\phi}=r_{\phi}$
and $I_{\phi,1}^{-1}=E[i_{\phi}^{2}(Z)].$ The efficiency bound is then
$E[i_{\phi}^{2}(Z)].$ It is known that in many cases bounds on $\left\Vert
r_{\phi}\right\Vert _{\beta}$ correspond to imposing smoothing conditions on
$r_{\phi}$ (see Kress 1999, Chapter 8). Hence, in these cases one can index
the level of irregularity by the level of smoothness of the influence function
$r_{\phi}.$ Related conditions have been extensively used in the literature of
ill-posed inverse problems in statistics and econometrics as \textquotedblleft
source conditions\textquotedblright, see e.g. Carrasco, Florens and Renault
(2007) and Chen and Reiss (2011). The novelty of Theorem \ref{ThmSVD} is in
the relation between \textquotedblleft smoothness\textquotedblright\ of
$r_{\phi}$ and identification and information.

\section{Semiparametric Models}

\label{SemiparametricModels}

This section studies the important class of semiparametric models, where
$\mathcal{P}=\{\mathbb{P}_{\theta,\eta}:\theta\in\Theta,\eta\in H\}.$ The
parameter space $\Lambda=\{(\theta,\eta):\theta\in\Theta,\eta\in H\}$ is a
subset of a Hilbert space $\mathbf{H}=\mathbb{R}^{p}\times\mathcal{H}.$ Define
$\langle(\theta_{1},\eta_{1}),(\theta_{2},\eta_{2})\rangle_{\mathbf{H}%
}:=\theta_{1}^{\prime}\theta_{2}+\langle\eta_{1},\eta_{2}\rangle_{\mathcal{H}%
}$. For semiparametric models the score operator, defined as in (\ref{3}), has
the representation (by the chain rule)%
\begin{equation}
S\left(  b_{\theta},b_{\eta}\right)  =\dot{l}_{\theta}^{\prime}b_{\theta}%
+\dot{l}_{\eta}b_{\eta},\text{ }b=(b_{\theta},b_{\eta})\in T(\lambda
_{0})\subseteq\mathbf{H}, \label{scoresemi}%
\end{equation}
where $\dot{l}_{\theta}\in L_{2}^{p}$ is the ordinary score function of
$\theta$ and $\dot{l}_{\eta}$ is a continuous linear operator from $T(\eta
_{0})\subset\mathcal{H}$ to $L_{2}.$ Let $\tilde{l}_{\theta}:=\dot{l}_{\theta
}-\Pi_{\overline{\mathcal{R}(\dot{l}_{\eta})}}\dot{l}_{\theta}$ be the
so-called efficient score function for $\theta$, where $\Pi_{\overline{V}}$
denotes the orthogonal projection operator onto $\overline{V}.$ The efficient
Fisher information matrix for $\theta$ is $\tilde{I}_{\theta}:=\mathbb{E}%
\left[  \tilde{l}_{\theta}\tilde{l}_{\theta}^{\prime}\right]  $. The following
result provides a characterization of the main identification condition for
the finite-dimensional parameter of a semiparametric model. For simplicity, I
consider the case $p=1,$ the extension to $p>1$ follows from applying the
result to the functionals $\phi(\lambda)=\alpha^{\prime}\theta$ for $\alpha
\in\mathbb{R}^{p}$. Section \ref{Nuisance} in the Supplemental Appendix
provides a parallel result for linear continuous functionals of the nuisance
parameter, allowing for $\theta$ to be infinite-dimensional and possibly
unidentified. Next proposition appears to be a new characterization of the
main condition for local identification in semiparametric models. \bigskip

\begin{proposition}
\label{Charactpar}For the functional $\phi(\lambda)=\theta\in\mathbb{R}$:
$\mathcal{N}(I_{\lambda_{0}})\subset\mathcal{N}(\dot{\phi})$ holds iff (i)
$\dot{l}_{\theta}\notin\overline{\mathcal{R}(\dot{l}_{\eta})}$ (positive
information $\tilde{I}_{\theta}>0)$ or (ii) $\dot{l}_{\theta}\in
\overline{\mathcal{R}(\dot{l}_{\eta})}\diagdown\mathcal{R}(\dot{l}_{\eta})$
(zero information $\tilde{I}_{\theta}=0$).
\end{proposition}

In the remainder of this section I extend some of the previous results to a
class of models that are nonlinear in the parameter of interest but linear in
nuisance parameters. \bigskip

\noindent\textbf{Assumption 2}: For some $\rho\geq1$ and for all
$\varepsilon>0,$ there exists $\delta>0$ and a continuous linear operator $S$
such that$,$ $\left\Vert \left(  f_{\lambda}-f_{\lambda_{0}}\right)
/f_{\lambda_{0}}-S(\lambda-\lambda_{0})\right\Vert <\varepsilon\left\vert
\theta-\theta_{0}\right\vert ^{\rho},$ for all $\lambda=(\theta,\eta
)\in\mathcal{B}_{\delta}(\lambda_{0}).$\bigskip

\noindent Assumption 2 is a mean-square differentiability condition with a
Lipschitz property on the derivative. It generally holds for models that are
nonlinear and smooth in the parameter of interest $\theta,$ but linear in the
nuisance parameters. Examples of models satisfying Assumption 2 include, among
others, structural models of unemployment duration in Heckman and Singer
(1984a, 1984b); linear and nonlinear panel data models with fixed effects (see
e.g. Bonhomme 2012); incomplete and complete games with multiple equilibria
(see e.g. Bajari, Hahn, Hong and Ridder 2011); semiparametric measurement
error models (see e.g. Hu and Schennach 2008); dynamic discrete choice models
(see e.g. Hu and Shum 2012); and binary discrete choice models with single and
multiple agents (see e.g. Chamberlain 1986, and more recently, Khan and
Nekipelov 2018). Importantly, Assumption 2 allows the nonparametric parameter
$\lambda_{0}$ to be unidentified and the parameter $\theta_{0}$ to be locally
irregularly identified, as it occurs in many of the aforementioned
applications. The latter feature differentiates our analysis from Chen et
al.'s (2014) setting. In most cases $1\leq\rho\leq2$ in Assumption 2, which
will limit the degree of irregularity permitted in identifying $\theta_{0},$
but functionals of the nuisance parameter are allowed to have arbitrary
degrees of irregularity, which can be important to accommodate many economic
applications with smooth densities (e.g. Heckman and Singer 1984a, 1984b).

The generalized Fisher information for $\theta$ is%

\[
I_{\theta,\rho}=\inf_{b\in\mathcal{B}_{\theta}}\frac{||Sb||^{2}}{\left\vert
\theta-\theta_{0}\right\vert ^{2\rho}},
\]
where $\mathcal{B}_{\theta}:=\{b\in\overline{T(\lambda_{0})}:b=(\theta
-\theta_{0},b_{\eta}),\theta\neq\theta_{0},$ $\left\vert \theta-\theta
_{0}\right\vert \leq1\}.$ It is straightforward to show that $\tilde
{I}_{\theta}=I_{\theta,1}$. Next theorem extends Theorem \ref{TheoremR} to the
nonlinear setting of Assumption 2.

\begin{theorem}
\label{NonlinearPositive}Let Assumption 2 hold. If $I_{\phi,\rho}>0$ for some
$\rho,$ $1\leq\rho<\infty,$ then $\theta_{0}\ $is locally identified:
regularly if $\rho=1$ and irregularly if $\rho>1$ and $\tilde{I}_{\theta}=0$.
\end{theorem}

\noindent Theorem \ref{NonlinearPositive} extends Theorem 7 in Chen et al.
(2014) to the semiparametric irregular case $\rho>1$. If the generalized
information $I_{\theta,\rho}$ is zero, Theorem \ref{Impossibility} implies
impossibility results on rates of convergence. Assumption 2 facilitates the
verification of the conditions for Theorem \ref{Impossibility} to hold. To see
this, consider a path $\lambda_{t}\in\Lambda$ passing through $\lambda_{0}$
such that $\left\Vert S\left(  \lambda_{t}-\lambda_{0}\right)  \right\Vert
=o(t^{\rho})$ and $\left\vert \theta_{t}-\theta_{0}\right\vert =Ct.$ For such
a path, Assumption 2 yields the conditions of Theorem \ref{Impossibility}.
Nevertheless, Assumption 2 is not necessary for Theorem \ref{Impossibility} to hold.

\section{Examples}

\label{Examplesdetail}

\subsection{Unemployment Duration with Nonparametric Heterogeneity}

\label{ExampleUnemp}

Nonparametric heterogeneity has played a critical role in rationalizing
unemployment duration ever since the seminal contributions by Heckman and
Singer (1984a, 1984b). Recent work by Alvarez et al. (2016) is motivated from
this perspective. These authors have shown nonparametric identification of the
distribution of unobserved heterogeneity, denoted by $G_{0},$ in their
nonparametric structural model for unemployment with two-spells under a sign
restriction on the parameter $\alpha,$ either $\alpha\geq0$ or $\alpha\leq0$.
As discussed by these authors, assuming either case imposes unattractive
restrictions on the economic model. For example, $\alpha\geq0$ implies that
all workers return to work eventually. With this background in mind, I
characterize the nonparametric identified set without prior restrictions. The
characterization of the identified set can be used to investigate alternative,
more attractive, conditions for nonparametric identification or to engage in
sensitivity identification analysis or partial identification bounds.

Let $\Lambda$ denote the parameter space for $\lambda,$ which is a subset of
densities wrt $\pi$. By Remark 1(i), it is natural to consider the score
operator defined on $T(\lambda_{0})\subset L_{1}(\pi),$ where $T(\lambda_{0})$
is the linear span of elements in $B_{0}=\{b\in L_{1}^{0}(\pi):\lambda
_{0}+b\in\Lambda\}.$ The score operator is, up to an irrelevant constant,
given by%
\[
Sb=\frac{1}{f_{\lambda_{0}}(t_{1},t_{2})}\int\frac{\beta^{2}}{t_{1}^{3/2}%
t_{2}^{3/2}}e^{-\frac{\left(  \alpha t_{1}-\beta\right)  ^{2}}{2t_{1}}%
-\frac{\left(  \alpha t_{2}-\beta\right)  ^{2}}{2t_{2}}}b(\alpha,\beta
)d\pi(\alpha,\beta).
\]
The following proposition, which builds on the nonparametric identification
results with sign restrictions in Alvarez et al. (2016), characterizes
$\mathcal{N}(S)$ under the following mild assumption.\bigskip

\noindent\textbf{Assumption 3}: (i) The set $\mathcal{T}\subseteq
\lbrack0,\infty)$ be a convex set with a non-empty interior; (ii) the measure
$\pi$ is such that $d\pi(-\alpha,\beta)=-d\pi(\alpha,\beta)$ and has no atom
at $\alpha=0$.\bigskip

\begin{proposition}
\label{Ex1} Under Assumption 3, $\mathcal{N}(S)=\left\{  b\in T(\lambda
_{0})\subset L_{1}^{0}(\pi):b(\alpha,\beta)=e^{-4\alpha\beta}b(-\alpha
,\beta)\right\}  .$
\end{proposition}

\noindent A corollary of this proposition is that the identified set for
$\lambda_{0}$ is the set $\{\lambda_{0}+b,$ where $b\in\mathcal{N}(S)\cap
B_{0}\}.$ An equivalent characterization of $\mathcal{N}(S)$ is
\[
\mathcal{N}(S)=\left\{  b\in T(\lambda_{0}):b(\alpha,\beta)=\frac
{C(\alpha,\beta)}{1-e^{4\alpha\beta}},\text{ where }C(\alpha,\beta)\text{ is
an odd function of }\alpha\right\}  .
\]
Proposition \ref{Identification} can then be used to check if a given linear
functional is identified or not. The characterization of $\mathcal{N}(S)$ also
can be used to find new point-identification results. The following corollary,
which follows directly from Proposition \ref{Ex1} illustrates this point.

\begin{corollary}
Under Assumption 3, if $\Lambda$ is a set of symmetric densities in $\alpha,$
i.e. for all $\lambda\in\Lambda,$ $\lambda(\alpha,\beta)=\lambda(-\alpha
,\beta)$ $\pi-$a.s, then $\lambda_{0}$ is nonparametrically identified.
\end{corollary}

The analysis of $\mathcal{N}(S)$ does not reveal, however, the degree of
identification for a given functional, whether regular or irregular. To
understand this, one must analyze the adjoint score operator. For this, it is
convenient to consider the approach of Remark 1(ii), which, as shown below,
leads to the score adjoint operator%
\begin{equation}
S^{\ast}g=\mathbb{E}\left[  \left.  g(Z)\right\vert \alpha,\beta\right]
,\text{ }g\in L_{2}^{0}.\label{adjoint}%
\end{equation}
Recall that under Remark 1(ii) $T(\lambda_{0})$ is the linear span of
$B_{0}=\{b\in L_{2}^{0}(G_{0}):\lambda_{0}+\lambda_{0}b\in\Lambda\}$, and one
say heterogeneity is nonparametric if $T(\lambda_{0})$ is dense in $L_{2}%
^{0}(G_{0})$ (i.e. $\overline{T(\lambda_{0})}=L_{2}^{0}(G_{0})$). The
following result provides a necessary condition for regular identification of
continuous linear functionals in this example.\bigskip

\begin{proposition}
\label{PropABS}Let Assumption 3(i) hold. If heterogeneity is nonparametric,
then $S^{\ast}$ has the representation in (\ref{adjoint}) and%
\[
\mathcal{R}(S^{\ast})\subset\left\{  b(\alpha,\beta)\in L_{2}^{0}%
(G_{0}):b(\alpha,\beta)=C_{1}+C_{2}\beta^{2}e^{2\alpha\beta}h(\alpha^{2}%
,\beta^{2})\right\}  ,
\]
for constants $C_{1}$ and $C_{2}$ and a continuous function $h(u,v)$ defined
on $(0,\infty)^{2}$ that, if $\mathcal{T}$ is bounded, is an infinite number
of times differentiable at $u\in(0,\infty),$ for all $v\in(0,\infty).$\bigskip
\end{proposition}

\noindent This result has important implications on which functionals are
regularly identified in this model. Take, for example, a moment of unobserved
heterogeneity%
\[
\phi(\lambda_{0})=\mathbb{E}\left[  r(\alpha,\beta)\right]  ,
\]
where $r\in L_{2}(G_{0}).$ From the reparametrization of Remark 1(ii), with
$\mathbf{H}=L_{2}(G_{0}),$
\[
\dot{\phi}(b)=%
%TCIMACRO{\dint }%
%BeginExpansion
{\displaystyle\int}
%EndExpansion
r(\alpha,\beta)b(\alpha,\beta)\lambda_{0}(\alpha,\beta)d\pi(\alpha
,\beta)=\langle r,b\rangle_{\mathbf{H}},
\]
and hence $r_{\phi}=r-\phi(\lambda_{0})\in L_{2}^{0}(G_{0}).$ Thus, only
smooth moment functions $r$ satisfying the conclusions of Proposition
\ref{Identification} will be regularly identified. Alvarez et al.'s (2016)
nonparametric identification result must imply $\mathcal{N}(S)=\{0\}$. In
turn, this result and Theorem 3 in Luenberger (1997, p.157) yield that
$\overline{\mathcal{R}(S^{\ast})}=L_{2}^{0}(G_{0}).$ Thus, the set
$\overline{\mathcal{R}(S^{\ast})}\setminus\mathcal{R}(S^{\ast})=L_{2}%
^{0}(G_{0})\setminus\mathcal{R}(S^{\ast})$ is rather large by Proposition
\ref{PropABS}, which shows that the class of irregularly identified
functionals is large in this model. Intuitively, this follows because the
density $f_{z/\alpha,\beta}(t_{1},t_{2})$ is very smooth in the parameters
$(\alpha,\beta)$, so that $\mathcal{R}(S^{\ast})$ only contains very smooth
functions. A direct implication of Proposition \ref{PropABS} is that the
cumulative distribution function (cdf) of unobserved heterogeneity at the
fixed point $(\alpha_{0},\beta_{0}),$ i.e. $\phi(\lambda_{0})=\mathbb{E}%
\left[  1(\alpha\leq\alpha_{0})1(\beta\leq\beta_{0})\right]  ,$ is not
regularly identified because $1(\alpha\leq\alpha_{0})1(\beta\leq\beta_{0})$ is
not continuous when $(\alpha_{0},\beta_{0})$ is in the interior of the support
of $\lambda_{0}$. Recall, $\phi(\lambda_{0})=\mathbb{E}\left[  1\left(
\alpha\leq\alpha_{0}\right)  1\left(  \beta\leq\beta_{0}\right)  \right]  $
for a fixed $\alpha_{0}<0<\beta_{0}\ $is the proportion of individuals at risk
of severe long term unemployment, and this proportion is thus non-regularly
identified, as an implication of Proposition \ref{PropABS}. We formalize the
discussion in the following corollary.

\begin{corollary}
\label{PropABScor}Let Assumption 3(i) hold. If heterogeneity is nonparametric
and $(\alpha_{0},\beta_{0})$ is in the interior of the support of $\lambda
_{0},$ then $\phi(\lambda_{0})=\mathbb{E}\left[  1\left(  \alpha\leq\alpha
_{0}\right)  1\left(  \beta\leq\beta_{0}\right)  \right]  $ is not regularly identified.
\end{corollary}

\subsection{Willingness-to-Pay}

\label{WTP}

The observed data is $Z_{i}=(Y_{i},V_{i},X_{i}),$ where $Y_{i}=1\left(
W_{i}\leq V_{i}\right)  ,$ $W_{i}$ is the unobserved willingness-to-pay of
individual $i,$ $V_{i}$ is a continuous observed random variable with known
and absolutely continuous conditional distribution $F_{V/X}$, and $X_{i}$ a
$d-$dimensional vector of covariates. The support of $W$ is $\mathcal{S}%
_{W}:=[0,w_{\max}]$ and that of $V$ is $\mathcal{S}_{V}:=[0,v_{\max}],$ with
$0<w_{\max},v_{\max}\leq\infty.$ It is assumed that $W_{i}$ and $V_{i}$ are
conditional independent given $X_{i}.$ Assume $G_{0}\left(  v,x\right)
=\mathbb{P}\left[  \left.  W_{i}\leq v\right\vert X_{i}=x\right]  $ is
differentiable at $v$ for each $x,$ and define $\lambda_{0}\left(  v,x\right)
=\partial G_{0}\left(  v,x\right)  /\partial v.$ The conditional density
$\lambda_{0}$ belongs to the parameter space $\Lambda,$ a subset of Lebesgue
densities on $\mathcal{S}_{W}\times\mathcal{S}_{X}.$ Define the measure
$\mu\left(  \left\{  0\right\}  \times B\right)  =\mu\left(  \left\{
1\right\}  \times B\right)  =\mu_{VX}(B),$ where $B$ is a Borel set of
$\mathbb{R}^{d+1}$ and $\mu_{VX}(B)$ is the probability measure for $(V,X).$
The density of $Z_{i}$ wrt $\mu$ is%
\begin{equation}
f_{\lambda_{0}}(y,v,x)=\left[  G_{0}\left(  v,x\right)  \right]  ^{y}\left[
1-G_{0}\left(  v,x\right)  \right]  ^{1-y}.\label{00}%
\end{equation}
Suppose one is interested in identification of the moment functional
$\phi(\lambda_{0})=\mathbb{E}\left[  r(W_{i},X)\right]  ,$ as in Lewbel,
McFadden and Linton (2011). I consider first this case, and then analyze
identification of the median WTP$,$%
\[
\phi(\lambda_{0})=Median(W)=\inf\{w\in\mathcal{S}_{W}:G_{0}\left(  w\right)
\geq1/2\},
\]
where $G_{0}\left(  w\right)  =\mathbb{P}\left[  W_{i}\leq w\right]  $ is the
cdf of $W_{i}.$

Following Remark 1(ii), the score operator is defined on $L_{2}^{0}(\pi),$
where $\pi$ denotes the probability measure for $(W,X)$ with conditional cdf
$G_{0}\left(  v,x\right)  .$ From (\ref{00}) the score operator is given for
$y\in\{0,1\},$ $v\in\mathcal{S}_{V}$ and $x\in\mathcal{S}_{X}$ by
\begin{equation}
S(b)=\frac{1}{f_{\lambda_{0}}(z)}1(f_{\lambda_{0}}>0)\left[  2y-1\right]
\int_{0}^{v}b(w,x)\lambda_{0}\left(  w,x\right)  dw.\label{SWTP}%
\end{equation}
By the Fundamental Theorem of Calculus $\mathcal{N}(S)=\{b\in T(\lambda
_{0})\subset L_{2}^{0}(\pi):b(w,x)=0$ $\pi-$a.s. on $0\leq w\leq v_{\max}\}.$
If $v_{\max}\geq w_{\max}$ then $\mathcal{N}(S)=\{0\}$ and all linear
functionals of $\lambda_{0}$ are identified by Proposition
\ref{Identification}.

If $v_{\max}<w_{\max}$ identification is more complicated and depends on the
functional. For the functional $\phi(\lambda_{0})=\mathbb{E}\left[
r(W_{i},X)\right]  ,$%
\begin{align*}
\dot{\phi}(b)  &  :=\phi(\lambda_{0}+b\lambda_{0})-\phi(\lambda_{0})\\
&  =\int_{\mathcal{S}_{W}\times\mathcal{S}_{X}}r(w,x)b(w,x)\lambda_{0}\left(
w,x\right)  d\mu_{X}dw\\
&  =\langle r,b\rangle_{\mathbf{H}},
\end{align*}
where $\mathbf{H}=L_{2}(\pi)$ and $\mu_{X}$ is the measure of $X.$ The unique
Riesz's representer of $\dot{\phi}$ in $\overline{T(\lambda_{0})}$ is given by%
\[
r_{\phi}=\Pi_{\overline{T(\lambda_{0})}}r.
\]
A sufficient condition for identification (cf. Proposition
\ref{Identification}) is%
\[
\int_{v_{\max}}^{w_{\max}}\int_{\mathcal{S}_{X}}r_{\phi}(w,x)b(w,x)\lambda
_{0}\left(  w,x\right)  d\mu_{X}dw=0
\]
for all $b\in\overline{T(\lambda_{0})}\subset L_{2}^{0}(\pi).$ Taking
$b(w,x)=r_{\phi}(w,x)$ leads to the following result:\bigskip

\begin{proposition}
\label{WTP0}If $v_{\max}<w_{\max},$ a sufficient condition for identification
of $\phi(\lambda_{0})=\mathbb{E}\left[  r(W_{i},X)\right]  $ with $r\in
L_{2}(\pi)$ is%
\[
r_{\phi}(w,x)=0\text{ }\pi-\text{a.s on }v_{\max}<w\leq w_{\max}.
\]

\end{proposition}

\noindent Consider henceforth the general case of a nonparametric model for
WTP, so $\overline{T(\lambda_{0})}=L_{2}^{0}(\pi).$ In this case%
\[
r_{\phi}=r-\phi(\lambda_{0}).
\]
To investigate the degree of identification, whether irregular or regular, let
us compute the adjoint score operator. First, note by conditional independence%
\[
S(b)=\frac{1}{f_{\lambda_{0}}(z)}1(f_{\lambda_{0}}>0)\left[  2y-1\right]
\mathbb{E}\left[  \left.  1(W\leq V)b(W,X)\right\vert V=v,X=x\right]  .
\]
Considering separately the case $y=1$ and $y=0,$ we see that
\[
S(b)=\mathbb{E}\left[  \left.  b(W,X)\right\vert Y=y,V=v,X=x\right]
1(f_{\lambda_{0}}>0).
\]
It is then well-known that the adjoint operator of a conditional mean operator
is
\[
S^{\ast}g(w,x)=\mathbb{E}\left[  \left.  g(Y,V,X)\right\vert W=w,X=x\right]
,
\]
which after some simple algebra gives%
\[
S^{\ast}g(w,x)=\int_{0}^{v_{\max}}1(w\leq v)g(1,v,x)f_{V/X=x}(v)dv+\int%
1(w>v)g(0,v,x)f_{V/X=x}(v)dv,
\]
where $f_{V/X=x}(v)$ is the known conditional Lebesgue density of $V$ given
$X=x.$ A direct consequence of this representation of $S^{\ast}$ and the
Fundamental Theorem of Calculus is that $S^{\ast}g(w,x)$ is absolutely
continuous in $w,$ for each $x$. This result can be used to characterize when
moment functionals such as $\phi(\lambda_{0})=\mathbb{E}\left[  r(W_{i}%
,X)\right]  $ are regularly identified. The regularity condition $r=S^{\ast
}g,$ for some $g\in L_{2},$ implies by the Fundamental Theorem of Calculus%
\[
\frac{\partial r(w,x)}{\partial w}=\left(  g(0,w,x)-g(1,w,x)\right)
f_{V/X=x}(w),
\]
which has a solution%
\begin{equation}
g(y,v,x)=\left[  1-2y\right]  \frac{\partial r(v,x)}{\partial v}\frac
{1}{f_{V/X=x}(v)}1(f_{V/X=x}(v)>0), \label{sol}%
\end{equation}
in $L_{2}$, provided
\begin{equation}
\int_{0}^{v_{\max}}\int\left[  \frac{\partial r(v,x)}{\partial v}\right]
^{2}\frac{1}{f_{V/X=x}(v)}d\mu_{X}(x)dv<\infty. \label{genreg}%
\end{equation}
Next result summarizes the findings on moments.

\begin{proposition}
\label{WTP1}Assume $\overline{T(\lambda_{0})}=L_{2}^{0}(\pi)$ and $v_{\max
}\geq w_{\max}.$ Let $r\in L_{2}(\pi)$ be a moment function such that $r(w,x)$
is differentiable at $w,$ for each $x$. Then, (\ref{genreg}) is necessary and
sufficient for regular identification of $\phi(\lambda_{0})=\mathbb{E}\left[
r(W_{i},X)\right]  .$
\end{proposition}

This result has interesting implications for known results in the literature.
First, the regular estimator that results from the moment representation based
on the solution (\ref{sol})
\begin{align*}
\phi(\lambda_{0})  &  =\mathbb{E}\left[  g(Y,V,X)\right] \\
&  =\mathbb{E}\left[  \left[  1-2Y\right]  \frac{\partial r(V,X)}{\partial
v}\frac{1}{f_{V/X}(V)}\right]
\end{align*}
is related to (but different from) the estimator proposed in Lewbel (1997)
without covariates, which is given by
\[
\mathbb{E}\left[  \left[  1(V\geq0)-Y\right]  \frac{\partial r(V)}{\partial
v}\frac{1}{f_{V}(V)}\right]  .
\]
The arguments above then show that the sufficient finite variance condition
derived in Lewbel (1997) for asymptotic normality of his estimators turns out
to be also necessary.

When applied to the mean of $W,$ i.e. $r(w,x)=w,$ which is a leading example
in Khan and Tamer (2003), the necessary and sufficient condition for regular
identification becomes simply%
\begin{equation}
\int_{0}^{v_{\max}}\int\frac{1}{f_{V/X=x}(v)}d\mu_{X}(x)dv<\infty. \label{reg}%
\end{equation}
If supports of $W$ and $V$ are unbounded, i.e. $w_{\max}=v_{\max}=\infty,$ so
$f_{V/X=x}(u)$ vanishes in the tails, the last condition does not hold, which
gives Khan and Tamer's (2010) irregularity result using a different method of
proof (they compute least favorable distributions and Fisher information).
Note that irregularity can also happen with bounded supports, depending on the
density $f_{V/X=x}(v),$ and that the condition is necessary (results which
were not considered in Khan and Tamer (2003) for the mean $r(w,x)=w$).

I now turn into the more difficult problem of semiparametric identification of
the median WTP, $\phi(\lambda_{0})=Median(W).$ This functional is nonlinear,
and hence the results of Theorem \ref{MainTheoremS} are not readily
applicable. For local identification, we require the following
condition:\bigskip

\noindent\textbf{Assumption 4}: (i) $v_{\max}>\phi(\lambda_{0})$ (ii)
$G_{0}\left(  v\right)  $ is continuous on $\mathcal{S}_{V}$ and
differentiable in a neighborhood of $v=\phi(\lambda_{0})$ with a positive
derivative $\mathbb{E}\left[  \lambda_{0}\left(  v,X_{i}\right)  \right]
.$\bigskip

\noindent Assumption 4(i) is a support condition that is weaker than that for
identification of moments of WTP. Assumption 4(ii) is typical in quantile
identification. Under Assumption 4, the median\ WTP has an influence function%
\[
r_{\phi}(w,x)=\frac{-\left\{  1(w<\phi(\lambda_{0}))-0.5\right\}  }%
{\mathbb{E}\left[  \lambda_{0}\left(  \phi(\lambda_{0}),X_{i}\right)  \right]
}.
\]
As the result above suggests, the discontinuity of the influence function
implies irregular identification. Next result, formalizes this finding.

\begin{proposition}
\label{WTP2}Let Assumption 4 hold, and assume $\overline{T(\lambda_{0})}%
=L_{2}^{0}(\pi).$ Then, $\phi(\lambda_{0})=Median(W)$ is irregularly
identified$.$
\end{proposition}

\subsection{Asset Pricing Euler Equation with Measurement Error}

\label{ExampleCAPM}

The goal of this example is to provide primitive conditions based on the
results of this paper for identification of the discount factor $\theta_{0}$
and measures of risk aversion. These are important parameters in these models.
For example, discount factors are a key determinant of individual's
intertemporal decisions such as asset accumulation (Venti and Wise, 1998,
Samwick, 2006), labor supply decisions (MaCurdy, 1981) and job search
(Dellavigna and Paserman, 2005).

\subsubsection{Identification of the Discount Factor}

The Euler equation with measurement error is a nonlinear conditional moment
restriction model. The first step in our analysis is to parametrize the model
in a way that makes it amenable to the results of this paper. To that end, I
consider the following assumption. Recall $C_{t}=m(C_{t}^{\ast},\varepsilon
_{t}).$\bigskip

\noindent\textbf{Assumption 5}: (i) $C_{t+1}^{\ast}$ is independent of
$C_{t},$ conditional on $C_{t+1}$, and measurable wrt $\sigma(C_{t+1}%
,\varepsilon_{t+1});$ (ii) $(\varepsilon_{t+1},\varepsilon_{t})$ is
independent of $R_{t+1},$ given $(C_{t+1},C_{t})$; (iii) the distribution of
$C_{t}^{\ast}$ conditional on $C_{t}$ does not depend on $t$.\bigskip

\noindent Assumption 5(i) can be relaxed to: $C_{t+1}^{\ast}$ is independent
of $C_{t-1},$ conditional on $(C_{t+1},C_{t}),$ at the cost of increasing the
dimension of the arguments in the nonparametric component given below. If $m$
is monotone in $C_{t}^{\ast},$ then 5(i) can be written in terms of
$\varepsilon_{t+1}.$ Assumption 5(i-iii) is less restrictive than typical
assumptions considered in the literature, which assume, in addition to
functional form assumptions on $m,$ that $\varepsilon_{t}$ is independent of
\textquotedblleft everything\textquotedblright; see, for example, Altonji and
Siow (1987), Runkle (1991), Dynam (2000), and Alan et al. (2009).

For the sake of exposition, I consider the case without household's
characteristics $X_{t,i}.$ The presence of $X_{t,i}$ in $\mathcal{F}_{t,i}$
adds additional moment restrictions, so it is simpler for identification. All
the arguments below can be easily adapted to the presence of $X_{t,i}.$ It is
also straightforward to extend the identification results to models with more
than one asset, habit formation or other observable variables in the marginal utility.

Assumption 5 ensures the following parametrization in terms of observables
$(C_{t+1},R_{t+1},C_{t}),$%
\begin{align*}
\mathbb{E}\left[  \left.  \theta_{0}\dot{u}(C_{t+1}^{\ast})R_{t+1}-\dot
{u}(C_{t}^{\ast})\right\vert C_{t}\right]   &  =\mathbb{E}\left[  \left.
\theta_{0}\dot{u}(C_{t+1}^{\ast})\mathbb{E}\left[  \left.  R_{t+1}\right\vert
C_{t+1},C_{t},\varepsilon_{t+1},\varepsilon_{t}\right]  -\dot{u}(C_{t}^{\ast
})\right\vert C_{t}\right] \\
&  =\mathbb{E}\left[  \left.  \theta_{0}\dot{u}(C_{t+1}^{\ast})\mathbb{E}%
\left[  \left.  R_{t+1}\right\vert C_{t+1},C_{t}\right]  -\dot{u}(C_{t}^{\ast
})\right\vert C_{t}\right]  ,\\
&  =\mathbb{E}\left[  \left.  \theta_{0}\mathbb{E}\left[  \left.  \dot
{u}(C_{t+1}^{\ast})\right\vert C_{t+1}\right]  \mathbb{E}\left[  \left.
R_{t+1}\right\vert C_{t+1},C_{t}\right]  -\dot{u}(C_{t}^{\ast})\right\vert
C_{t}\right] \\
&  =\mathbb{E}\left[  \left.  \theta_{0}\eta_{0}(C_{t+1})R_{t+1}-\eta
_{0}(C_{t})\right\vert C_{t}\right]  ,
\end{align*}
where $\eta_{0}(C_{t+1})=\mathbb{E}\left[  \left.  \dot{u}(C_{t+1}^{\ast
})\right\vert C_{t+1}\right]  $ and the first equality uses 4(i), the second
4(ii), the third 4(i) and the fourth 4(iii). This new parametrization is a
nonlinear conditional moment indexed by the discount factor $\theta_{0}$ and
the projected marginal utility $\eta_{0}.$ Let $\mu$ be the probability
measure of $C_{t}.$ Then, I denote $\lambda_{0}=(\theta_{0},\eta_{0}%
)\in\Lambda=\Theta\times H,$ $\Theta\subset(0,1)\mathbb{\ }$and $H\subset
L_{2}(\mu)$.

The following condition guarantees that the conditional mean operator%
\begin{equation}
A\eta(c)=\mathbb{E}\left[  \left.  \eta(C_{t+1})R_{t+1}\right\vert
C_{t}=c\right]  , \label{Fredholm2}%
\end{equation}
is well-defined and compact when viewed as $A:L_{2}(\mu)\rightarrow L_{2}%
(\mu)$. This is a standard assumption in the literature, see e.g. Carrasco et
al. (2007). Let $g(C_{t+1},C_{t})$ be the joint Lebesgue density of
$(C_{t+1},C_{t}),$ and let $f_{t+1}(C_{t+1})$ and $f_{t}(C_{t})$ denote its
marginals, respectively.\bigskip

\noindent\textbf{Assumption 5}: (iv) $0<\mathbb{E}\left[  R_{t+1}^{2}\left\{
g(C_{t+1},C_{t})/f_{t}(C_{t+1})f_{t}(C_{t})\right\}  \right]  <\infty.$
\bigskip

\noindent Define, for $\lambda=(\theta,\eta)\in\Lambda=\Theta\times H,$%
\[
M_{\lambda}=\theta A\eta(c)-\eta(c).
\]
The identification results of Section \ref{SemiparametricModels} can be
adapted to this moment model replacing $\left(  f_{\lambda}-f_{\lambda_{0}%
}\right)  /f_{\lambda_{0}}$ by $M_{\lambda_{0}+b}-M_{\lambda_{0}}$ and the
norm $\left\Vert \cdot\right\Vert $ by $\left\Vert \cdot\right\Vert _{2,\mu}.$
I then proceed to verify Assumption 2 in this example. It is straightforward
to show that for all $b=(b_{\theta},b_{\eta})\in\mathbb{R}\times L_{2}(\mu),$
with $b_{\theta}=\theta-\theta_{0}$%
\begin{equation}
\left\Vert M_{\lambda_{0}+b}-M_{\lambda_{0}}-S(b)\right\Vert _{2,\mu}%
\leq\left\Vert \theta_{0}Ab_{\eta}(c)-b_{\eta}(c)\right\Vert _{2,\mu
}\left\vert b_{\theta}\right\vert , \label{A2dis}%
\end{equation}
where%
\[
S\left(  b_{\theta},b_{\eta}\right)  =\dot{l}_{\theta}b_{\theta}+\dot{l}%
_{\eta}b_{\eta},
\]%
\[
\dot{l}_{\theta}b_{\theta}=b_{\theta}A\eta_{0}\text{ and }\dot{l}_{\eta
}b_{\eta}=\theta_{0}Ab_{\eta}(c)-b_{\eta}(c).
\]
Since $A$ is bounded, for each $\varepsilon>0,$ one can make $\left\Vert
\theta_{0}Ab_{\eta}(c)-b_{\eta}(c)\right\Vert _{2,\mu}<\varepsilon$ by
choosing $\delta$ small enough and $\left\Vert b_{\eta}\right\Vert _{2,\mu
}<\delta.$ Thus, Assumption 2 holds with $\rho=1.$

It follows from the previous parametrization and Assumption 5 that local
identification of the discount factor is regular.\footnote{The actual
effective score operator of the model is proportional to $S,$ with a
proportionality \textquotedblleft constant\textquotedblright\ given by the
inverse of the conditional variance $\mathbb{E}\left[  \left.  V_{t+1}%
^{2}\right\vert C_{t}=c\right]  ,$ where $V_{t+1}=\theta_{0}\eta_{0}%
(C_{t+1})R_{t+1}-\eta_{0}(C_{t}).$ This conditional variance is assumed to be
bounded and bounded away from zero. Thus, statements related to regularity can
be given in terms of $S$.} Formally, Theorem 3.2 in Kress (1999, p. 29)
implies that the range of $\dot{l}_{\eta},$ $\mathcal{R}(\dot{l}_{\eta})$, is
closed. It follows from Rudin (1973, 4.14) that $\mathcal{R}(\dot{l}_{\eta
}^{\ast})$ is also closed, and by the expression above, $\mathcal{R}(S^{\ast
})$ is also closed. Then, by the results of this paper, all locally identified
linear continuous functionals of $\lambda_{0}=(\theta_{0},\eta_{0})$ are
regularly locally identified. Note that this does not apply to functionals of
the marginal utility $\dot{u}_{0}.$

Positive information for the discount factor $\theta_{0}$ holds iff
$(1,0)\in\mathcal{R}(S^{\ast})$, which means there exists $g\in L_{2}(\mu)$
such that%
\begin{equation}
\langle A\eta_{0},g\rangle=1,\dot{l}_{\eta}^{\ast}g=0. \label{RegCAPM}%
\end{equation}
Since $\dot{l}_{\eta}^{\ast}g=\theta_{0}A^{\ast}g(c)-g(c),$ the equation
$\dot{l}_{\eta}^{\ast}g=0$ means that $g$ is an eigenfunction of $A^{\ast}$
with eigenvalue $\theta_{0}^{-1}.$ Such eigenfunction exists because
eigenvalues of $A^{\ast}$ are complex conjugates of those of $A$ and
$\theta_{0}$ is real-valued. Then, a sufficient condition for local
identification of the discount factor is that for such eigenfunction, say
$g_{0}:$ $\langle\eta_{0},g_{0}\rangle\neq0$. Note that $\langle\eta_{0}%
,g_{0}\rangle\neq0$ and $\theta_{0}>0$ guarantee (\ref{RegCAPM}) by choosing
$g=cg_{0}$ with $c=1/\langle A\eta_{0},g_{0}\rangle$ (since multiples of
eigenfunctions are eigenfunctions). The discussion is summarized in the
following result.

\begin{proposition}
\label{discount}Let Assumption 5 hold and assume $\langle\eta_{0},g_{0}%
\rangle\neq0$ for $g_{0}$ one of the eigenfunctions of $A^{\ast}$
corresponding to the eigenvalue $\theta_{0}^{-1}.$ Then, $\theta_{0}$ is
locally regularly identified.
\end{proposition}

The condition $\langle\eta_{0},g_{0}\rangle\neq0$ is mild, and holds, for
example, when $\eta_{0}$ and $g_{0}\ $are positive. Escanciano and Hoderlein
(2010) present primitive conditions for nonparametric identification of
positive $\eta_{0}$ and $g_{0}\ $based on Perron-Frobenius theory. Chen et al
(2014) and Escanciano et al. (2015) also use Perron-Frobenius to obtain
identification of related but different Euler equation models. See also Hansen
and Scheinkman (2009) and Christensen (2017) for other applications of
Perron-Frobenius theory. These identification results are nonparametric and
for models without measurement error. In contrast, a simple
\textit{semiparametric} identification condition for the discount factor is
presented here, in a model with measurement error, namely $\langle\eta
_{0},g_{0}\rangle\neq0.$ Proposition \ref{discount} thus shows that regular
local identification of the discount factor holds under rather general
conditions on the measurement error mechanism ($m$ and the distribution of
$\varepsilon_{t}$ are nonparametric and unidentified under our conditions).

An important empirical literature has provided estimation and inference
results on Euler equations accounting for measurement error. Papers within
this literature use functional form assumptions for utilities and for the
measurement error mechanism. The identification result of Proposition
\ref{discount} opens the door for more robust empirical strategies for
inference on the discount factor in microeconomic applications based on the
Euler Equation. For example, Altonji and Siow (1987), Runkle (1991), Dynam
(2000), and Alan, Attanasio and Browning (2009) assume parametric marginal
utilities and $m(C_{t}^{\ast},\varepsilon_{t})=C_{t}^{\ast}\varepsilon_{t},$
i.e.%
\begin{equation}
C_{t}=C_{t}^{\ast}\varepsilon_{t}, \label{multiplicative}%
\end{equation}
with $\varepsilon_{t}$ independent of everything and, in some cases, assumed
to be log normally distributed. The identification result above shows that
regular local identification of the discount factor follows under more general
assumptions than previously recognized, including situations where the
marginal utility, the measurement equation and other nonparametric parameters
are not identified. This point illustrates the concept of semiparametric
identification emphasized in this paper.

\subsubsection{Identification of Average Risk Aversion}

The Average Arrow-Pratt coefficient of Absolute Risk Aversion (AARA) parameter
is given by%
\[
\chi(\dot{u}_{0})=\mathbb{E}\left[  \frac{-\partial\dot{u}_{0}(C_{t}^{\ast
})/\partial C_{t}^{\ast}}{\dot{u}_{0}(C_{t}^{\ast})}\right]  .
\]
The following conditions guarantee that this parameter is well-defined, and
satisfies some properties given below. Let $\mu^{\ast}$ denote the probability
measure of $C_{t}^{\ast}$, with density $f_{t}^{\ast}(\cdot),$ and assume the
parameter space for marginal utilities $\dot{U}$ satisfies, for a small
positive number $\epsilon,$
\[
\dot{U}\subset\left\{  \dot{u}\in L_{2}(\mu^{\ast}):\dot{u}(c^{\ast}%
)\geq\epsilon>0\text{ and }\int\dot{u}^{2}d\mu^{\ast}=1\right\}  .
\]
The scale normalization in $\dot{U}$ is necessary to achieve identification in
nonparametric Euler equations.\bigskip

\noindent\textbf{Assumption 5}: (v) The functions $\log\left(  \dot{u}%
_{0}(\cdot)\right)  $ and $\log(f_{t}^{\ast}(\cdot))$ are continuously
differentiable on the convex support of $C_{t}^{\ast}$ (possibly unbounded)
and these functions and their derivatives are in $L_{2}(\mu^{\ast}).$ The true
marginal utility satisfies $\dot{u}_{0}\in\dot{U}.$ \bigskip

\noindent By Assumption 5(v) and integration by parts,%
\[
\chi(\dot{u}_{0})=\mathbb{E}\left[  \log\left(  \dot{u}_{0}(C_{t}^{\ast
})\right)  d(C_{t}^{\ast})\right]  ,
\]
where
\begin{equation}
d(c^{\ast})\equiv\frac{\partial f_{t}^{\ast}(c^{\ast})}{\partial c^{\ast}%
}\frac{1}{f_{t}^{\ast}(c^{\ast})}. \label{scores}%
\end{equation}
The functional $\chi(\dot{u}_{0}),$ although nonlinear, is concave and
differentiable, with pathwise derivative
\[
\dot{\chi}(b)=\mathbb{E}\left[  b(C_{t}^{\ast})\frac{d(C_{t}^{\ast})}{\dot
{u}_{0}(C_{t}^{\ast})}\right]  .
\]
Thus, the AARA parameter as a functional on $L_{2}(\mu^{\ast})$ has the
Riesz's representer
\begin{equation}
r_{\chi}(C_{t}^{\ast})=\frac{d(C_{t}^{\ast})}{\dot{u}_{0}(C_{t}^{\ast})}.
\label{rAARA}%
\end{equation}
To link the marginal utility with the projected marginal utility we need to
define the conditional mean operator $L:L_{2}(\mu^{\ast})\rightarrow L_{2}%
(\mu),$%
\[
L\dot{u}(C_{t})=\mathbb{E}\left[  \left.  \dot{u}(C_{t}^{\ast})\right\vert
C_{t}\right]  ,
\]
which has an adjoint operator $L^{\ast}:L_{2}(\mu)\rightarrow L_{2}(\mu^{\ast
})$ given by%
\[
L^{\ast}w(C_{t})=\mathbb{E}\left[  \left.  w(C_{t})\right\vert C_{t}^{\ast
}\right]  .
\]
Applying the results of this paper, one obtains the following identification
result for the AARA.

\begin{proposition}
\label{AARA}Let Assumption 5 hold and assume $\theta_{0}$ is locally
identified. Then, the following condition is sufficient for regular local
identification of the AARA: there exists $g\in L_{2}(\mu)$ such that, with
$r_{\chi}$ given in (\ref{rAARA}) and $\dot{l}_{\eta}^{\ast}g=\theta
_{0}A^{\ast}g(c)-g(c)$,
\begin{equation}
r_{\chi}(C_{t}^{\ast})=\mathbb{E}\left[  \left.  \dot{l}_{\eta}^{\ast}%
g(C_{t})\right\vert C_{t}^{\ast}\right]  . \label{AARAreg}%
\end{equation}

\end{proposition}

\noindent It is convenient to decompose (\ref{AARAreg}) into two parts: (i)
existence of $w\in L_{2}(\mu)$ such that
\begin{equation}
r_{\chi}(C_{t}^{\ast})=\mathbb{E}\left[  \left.  w(C_{t})\right\vert
C_{t}^{\ast}\right]  ; \label{P1}%
\end{equation}
and (ii) conditions that guarantee that such $w$ belongs to $\mathcal{R}%
(\dot{l}_{\eta}^{\ast}).$ I provide simple primitive conditions for (\ref{P1})
to hold in the multiplicative measurement error model (\ref{multiplicative}),
when the error density is known and given by $f_{\varepsilon},$ e.g.
log-normal as in e.g. Alan, Attanasio and Browning (2009). In this model, the
regularity condition (\ref{P1}) is$,$
\[
r_{\chi}(c^{\ast})=\int f_{\varepsilon}(c/c^{\ast})w(c)dc.
\]
The following Lemma provides a sufficient condition for existence of $w\in
L_{2}(\mu)$ satisfying this equation. Let $L_{1}(\mathbb{R})$ and
$L_{2}(\mathbb{R})$ denote the set of integrable and squared integrable
functions, respectively. For $f\in L_{1}(\mathbb{R}),$ define the Fourier
transform $\hat{f}=(2\pi)^{-1/2}\int e^{-itx}f(t)dt,$ where $i=\sqrt{-1}.$
Define $K(u)=\exp(u)f_{\varepsilon}(\exp(u))$ and $x(\tau)=\exp(-\tau)r_{\chi
}(\exp(\tau)).$

\begin{lemma}
\label{LemmaP1}If $K(u)$ is symmetric in $u,$ $x\in L_{2}(\mathbb{R})$ and
$\hat{x}/\hat{K}\in L_{2}(\mathbb{R}),$ then, there exists a solution $w\in
L_{2}(\mu)$ of (\ref{P1}). Moreover, a solution is given by%
\[
w(c)=\frac{1}{2\pi}\operatorname{Re}\int e^{it\log(c)}\frac{\hat{x}(t)}%
{\hat{K}(t)}dt,
\]
where $\operatorname{Re}$ denotes the real part.
\end{lemma}

The symmetry condition on $K$ is satisfied by the log normal distribution used
in the empirical literature. I now provide primitive conditions for
$w\in\mathcal{R}(\dot{l}_{\eta}^{\ast})$. Note that if $w\notin\mathcal{R}%
(\dot{l}_{\eta}^{\ast})$ the AARA is not identified. By duality,
$w\in\mathcal{R}(\dot{l}_{\eta}^{\ast})=\mathcal{N}(\dot{l}_{\eta})^{\perp}$
has a simple interpretation: $w$ is orthogonal to all projected marginal
utilities solving the Euler equation, i.e.%
\begin{equation}
\mathbb{E}\left[  w(C_{t})\eta(C_{t})\right]  =0\text{ for all }\eta\text{
such that }\theta_{0}A\eta=\eta. \label{ortho}%
\end{equation}
By compactness of $A,$ the space of such $\eta^{\prime}s$ is
finite-dimensional (see Kress 1999), which means that (\ref{ortho}) can be
tested. Importantly, (\ref{ortho}) holds for $\eta_{0}$ under (\ref{P1}),
since by iterated expectations%
\begin{align*}
\mathbb{E}\left[  w(C_{t})\eta_{0}(C_{t})\right]   &  =\mathbb{E}\left[
r_{\chi}(C_{t}^{\ast})\eta_{0}(C_{t})\right] \\
&  =\mathbb{E}\left[  r_{\chi}(C_{t}^{\ast})\dot{u}_{0}(C_{t}^{\ast})\right]
\\
&  =\mathbb{E}\left[  d(C_{t}^{\ast})\right] \\
&  =0.
\end{align*}
A primitive condition for identification of $\eta_{0}$ is $r(C_{t+1},C_{t})>0$
and $g(C_{t+1},C_{t})>0,$ where $r(C_{t+1},C_{t})=\mathbb{E}\left[  \left.
R_{t+1}\right\vert C_{t+1},C_{t}\right]  .$ Thus, these primitive conditions
and the mild integrability conditions of Lemma \ref{LemmaP1} imply regular
identification of the ARRA by virtue of Proposition \ref{AARA}.

\section{Conclusions}

\label{Conclusions} This paper provides tools for investigating semiparametric
identification, with a particular emphasis on irregular identification. First,
it considers semiparametric identification for linear models and obtains
necessary and sufficient conditions for regular and irregular identification.
I then show that semiparametric irregular identification is a common feature
of many economic models of practical interest. The application to the
structural model for unemployment in Alvarez, Borovickov\'{a} and Shimer
(2016) illustrates this point. The Euler Equation application illustrates the
usefulness of the characterization of regular identification and its
applicability to the important class of conditional moment models. Regular
identification of the discount factor and measures of risk aversion can be
obtained under relatively simple conditions, despite the nonlinearity of both
the model and the AARA functional.

The question of whether zero information corresponds to a lack of
identification is a rather delicate question, as was first pointed out in
Chamberlain (1986). Indeed, I show here that irregularity corresponds
mathematically to a boundary case in an infinite-dimensional space. Regular
identification is, however, easier to characterize and, under mild smoothness
conditions, a positive semiparametric Fisher information for the parameter
implies its local identification. When the Fisher information is zero,
positivity of a new generalized Fisher information introduced in this paper
implies irregular identification. When the generalized Fisher information is
zero, I obtain impossibility results on rates of convergence. The
impossibility results on regular identification and rates apply to both linear
and nonlinear models and parameters.

A number of issues deserve further study. For example, it will be useful to
investigate primitive conditions for positive or zero generalized Fisher
information in specific economic applications and, using the tools provided
here, to see how these conditions translate into specific rates of convergence
for estimators. Likewise, Section \ref{Nonlinear} in the Supplemental Appendix
provides sufficient conditions for regular and irregular semiparametric
identification in nonlinear models and for nonlinear functionals. Applying
these results to specific examples, and establishing connections with
attainability of rates of convergence for semiparametric estimators, remain
topics for future research.

\section{Appendix: Proofs of Main Results}

\label{Proofs}

\noindent\textbf{Proof of Proposition \ref{Identification}}: I first show that
$\mathcal{N}(I_{\lambda_{0}})=\mathcal{N}(S).$ From the definition of
$I_{\lambda_{0}},$ the implication $\mathcal{N}(S)\subset\mathcal{N}%
(I_{\lambda_{0}})$ trivially holds. The other implication follows from
$\langle I_{\lambda_{0}}b,b\rangle_{\mathbf{H}}=\left\Vert Sb\right\Vert
^{2}.$ Having $\mathcal{N}(I_{\lambda_{0}})=\mathcal{N}(S),$ that (\ref{1}) is
a sufficient condition for identification of $\phi(\lambda_{0})$ follows from
the definition of identification. If $f_{\lambda}=f_{\lambda_{0}},$ then
$b=\lambda-\lambda_{0}\in\mathcal{N}(S)\subset\mathcal{N}(\dot{\phi}),$ and
hence $\phi(\lambda_{0}+b)=\phi(\lambda_{0}).$ To prove the necessity, suppose
that (\ref{1}) does not hold, i.e.%
\[
\mathcal{N}(S)\varsubsetneq\mathcal{N}(\dot{\phi}),
\]
then there exists $b\in T(\lambda_{0})$ such that $b\in\mathcal{N}(S)$ but
$b\notin\mathcal{N}(\dot{\phi}).$ This means by linearity that for all
$c\in\mathbb{R},$ $c\neq0,$ $cb\in\mathcal{N}(S)$ but $cb\notin\mathcal{N}%
(\dot{\phi}).$ By Assumption 1(iii) there exists $c$ such that $\lambda
_{0}+cb\in\Lambda$, $cb\in\mathcal{N}(S)$ and $cb\notin\mathcal{N}(\dot{\phi
}).$ This implies that $f_{\lambda_{0}+cb}=f_{\lambda_{0}}$ but $\phi
(\lambda_{0}+cb)\neq\phi(\lambda_{0}).$ That is, $\phi(\lambda_{0})$ is not
identified. $\blacksquare$\bigskip

\noindent\textbf{Proof of Theorem \ref{MainTheoremS}}: By Proposition 1 in
Luenberger (1997, p.52)%
\[
\mathcal{N}(S)\subset\mathcal{N}(\dot{\phi})
\]
is equivalent to
\[
\mathcal{N}(\dot{\phi})^{\perp}\subset\mathcal{N}(S)^{\perp},
\]
since both $\mathcal{N}(S)$ and $\mathcal{N}(\dot{\phi})$ are closed linear
subspaces, and where henceforth $V^{\perp}$ denotes the orthocomplement of the
subspace $V$. However, since
\[
\dot{\phi}(b)=\langle b,r_{\phi}\rangle_{\mathbf{H}},
\]
for all $b\in T(\lambda_{0}),$ it follows that $\mathcal{N}(\dot{\phi}%
)^{\perp}=span\{r_{\phi}\}.$ On the other hand, by Theorem 3 in Luenberger
(1997, p.157)%
\[
\mathcal{N}(S)^{\perp}=\overline{\mathcal{R}(S^{\ast})}.
\]
The identification part follows from Proposition \textbf{\ref{Identification}%
}. The qualification of regular or irregular follows from Theorem 4.1 in van
der Vaart (1991), which shows that $r_{\phi}\in\mathcal{R}(S^{\ast
})\Longleftrightarrow I_{\phi}>0.$ $\blacksquare$\bigskip

\noindent\textbf{Proof of Theorem \ref{TheoremR}}: Suppose $\phi(\lambda_{0})$
is not identified. That means we can find $\lambda\in\Lambda$ such that
$\phi(\lambda)\neq\phi(\lambda_{0})$ and $f_{\lambda}=f_{\lambda_{0}}.$ Then,
$b=\lambda-\lambda_{0}\in\mathcal{N}(S)$ but $\dot{\phi}(b)\neq0.$ Choose $c$
sufficiently small so that $\left\vert \dot{\phi}(cb)\right\vert \leq1$, and
hence $cb\in\mathcal{B}_{\phi}$ and%
\[
I_{\phi,\rho}=\inf_{b\in\mathcal{B}_{\phi}}\frac{||Sb||^{2}}{[\dot{\phi
}(b)]^{2\rho}}=0,
\]
contradicting the positivity of $I_{\phi,\rho}.$ $\blacksquare$\bigskip

\noindent\textbf{Proof of Theorem \ref{Impossibility}}: Fix $K>0$ and take
$t\equiv t_{n}=(2K/C)n^{-1/(2\rho)}$ and $\lambda_{n}=\lambda_{t_{n}},$ so
that $\left\vert \phi(\lambda_{n})-\phi(\lambda_{0})\right\vert \geq
2Kn^{-1/(2\rho)}$ and $\left\Vert \left(  f_{\lambda_{n}}-f_{\lambda_{0}%
}\right)  /f_{\lambda_{0}}\right\Vert ^{2}<\varepsilon(2K/C)^{2\rho}n^{-1}.$
Let $\mathbb{P}_{n}$ denote the probability measure of $f_{\lambda_{n}},$ and
$\mathbb{P}_{0}$ that of $f_{\lambda_{0}},$ and $\mathbb{P}_{n}^{n}$ and
$\mathbb{P}_{0}^{n}$ their corresponding $n-$fold product. Lemma 1 in Lecam
(1973) and the basic inequality $\left(  \sqrt{a}-\sqrt{b}\right)  ^{2}%
\leq(a-b)^{2}/b$ for $a,b>0,$ imply that for each $\delta>0$ we can find
$\gamma>0$ sufficiently small such that
\[
v(\mathbb{P}_{n}^{n},\mathbb{P}_{0}^{n})<\delta\text{ if }\left\Vert \left(
f_{\lambda_{n}}-f_{\lambda_{0}}\right)  /f_{\lambda_{0}}\right\Vert
^{2}<\gamma/n,
\]
where $v(\mathbb{P},\mathbb{Q})$ denotes the total variation distance between
$\mathbb{P}$ and $\mathbb{Q}.$ The details are given as follows. Lemma 1 in
Lecam (1973) gives the inequality%
\[
v(\mathbb{P}_{n}^{n},\mathbb{P}_{0}^{n})\leq y(2-y^{2})^{1/2}%
\]
if $H(\mathbb{P}_{n},\mathbb{P}_{0})\leq y/\sqrt{n}$ and $y\leq1.$ From the
basic inequality $\left(  \sqrt{a}-\sqrt{b}\right)  ^{2}\leq(a-b)^{2}/b$ for
$a,b>0,$ we have%
\begin{align*}
H^{2}(\mathbb{P}_{n},\mathbb{P}_{0})  &  =\frac{1}{2}\int\left(
\sqrt{f_{\lambda_{n}}}-\sqrt{f_{\lambda_{0}}}\right)  ^{2}d\mu\\
&  \leq\frac{1}{2}\int\left[  \left(  f_{\lambda_{n}}-f_{\lambda_{0}}\right)
^{2}/f_{\lambda_{0}}\right]  d\mu\\
&  =\frac{1}{2}\left\Vert \left(  f_{\lambda_{n}}-f_{\lambda_{0}}\right)
/f_{\lambda_{0}}\right\Vert ^{2}\\
&  <y^{2}/n,
\end{align*}
with $y^{2}=\varepsilon(2K/C)^{2\rho}/2$ (where the last inequality follows
from the assumptions of the theorem). Since we can always choose
$\varepsilon>0$ small enough so that $y\leq1,$ Lemma 1 in Lecam (1973) can be
applied to achieve $v(\mathbb{P}_{n}^{n},\mathbb{P}_{0}^{n})<\delta$ by making
$y$ small enough. Next, Lemma 8 in Ishwaran (1996) implies that the rate of
any estimator cannot be better than $O_{P}(2Kn^{-1/(2\rho)})$. Since this
holds for each fixed $K>0,$ it follows that the rate of any estimator must be
slower than $n^{-1/(2\rho)}.$ $\blacksquare$\bigskip

\noindent\textbf{Proof of Theorem \ref{ThmSVD}}: The first part follows from
Theorem \ref{MainTheoremS}(i) and the characterization of $\mathcal{R}%
(S^{\ast})$ as those elements $b\in\overline{T(\lambda_{0})}$ with $\left\Vert
b\right\Vert _{1}<\infty.$ Any element $b\in\overline{T(\lambda_{0})}$ has the
singular value expansion (cf. Kress, 1999, Theorem 15.16)
\[
b=\sum_{j=1}^{\infty}\left\langle b,\varphi_{j}\right\rangle _{_{\mathbf{H}}%
}\varphi_{j}+\Pi_{\mathcal{N}(S)}b,
\]
which implies under semiparametric identification of the functional%
\[
\dot{\phi}(b)=\sum_{j=1}^{\infty}\left\langle b,\varphi_{j}\right\rangle
_{_{\mathbf{H}}}\left\langle r_{\phi},\varphi_{j}\right\rangle _{_{\mathbf{H}%
}}%
\]
and%
\[
Sb=\sum_{j=1}^{\infty}\lambda_{j}\left\langle b,\varphi_{j}\right\rangle
_{_{\mathbf{H}}}\psi_{j}.
\]
By Cauchy-Schwarz, for $b\in T(\lambda_{0}),$
\begin{align*}
\left\vert \dot{\phi}(b)\right\vert  &  \leq\left(  \sum_{j=1}^{\infty}%
\lambda_{j}^{-2}\left\langle r_{\phi},\varphi_{j}\right\rangle _{_{\mathbf{H}%
}}^{2}\right)  ^{1/2}\left(  \sum_{j=1}^{\infty}\lambda_{j}^{2}\left\langle
b,\varphi_{j}\right\rangle _{_{\mathbf{H}}}^{2}\right)  ^{1/2}\\
&  =\left\Vert r_{\phi}\right\Vert _{1}||Sb||.
\end{align*}
This yields%
\[
I_{\phi,1}\geq\left\Vert r_{\phi}\right\Vert _{1}^{-2}.
\]
By regular identification there exists $i_{\phi}$ with finite variance such
that $S^{\ast}i_{\phi}=r_{\phi},$ and hence another version of Cauchy-Schwarz
gives%
\begin{align*}
\left\vert \dot{\phi}(b)\right\vert  &  =\left\langle Sb,i_{\phi}\right\rangle
_{_{\mathbf{H}}}\\
&  \leq\left\Vert i_{\phi}\right\Vert ||Sb||.
\end{align*}
Indeed, it can be shown that $\left\Vert i_{\phi}\right\Vert =\left\Vert
r_{\phi}\right\Vert _{1}$ and equality in Cauchy-Schwarz holds when $Sb$ and
$i_{\phi}$ are proportional. The formal argument requires constructing
\[
b_{k}=\frac{\sum_{j=1}^{k}\lambda_{j}^{-1}\left\langle i_{\phi},\psi
_{j}\right\rangle \varphi_{j}}{\sum_{j=1}^{k}\left\langle i_{\phi},\psi
_{j}\right\rangle ^{2}}.
\]
Note $\dot{\phi}(\varphi_{j})=\left\langle r_{\phi},\varphi_{j}\right\rangle
_{_{\mathbf{H}}}=\lambda_{j}\left\langle i_{\phi},\psi_{j}\right\rangle $
(since $S\varphi_{j}=\lambda_{j}\psi_{j})$ and hence%
\[
\left\vert \dot{\phi}(b_{k})\right\vert =1,
\]
and $b_{k}$ is in the set $\mathcal{B}_{\phi}=\left\{  b\in T(\lambda
_{0}):\dot{\phi}(b)\neq0,\left\vert \dot{\phi}(b)\right\vert \leq1\right\}  .$
Moreover,%
\[
||Sb_{k}||^{2}=\frac{1}{\sum_{j=1}^{k}\left\langle i_{\phi},\psi
_{j}\right\rangle ^{2}}%
\]
which can get arbitrarily close to the bound $\left\Vert i_{\phi}\right\Vert
^{-2}=\left\Vert r_{\phi}\right\Vert _{1}^{-2}$ by choosing a large $k.$ Thus,
we conclude $I_{\phi,1}=\left\Vert r_{\phi}\right\Vert _{1}^{-2}=\left\Vert
i_{\phi}\right\Vert ^{-2}.$

For the second part, note that $\left\Vert r_{\phi}\right\Vert _{\beta}%
<\infty$ implies $r_{\phi}\in\overline{\mathcal{R}(S^{\ast})},$ and hence
identification. The condition $\left\Vert r_{\phi}\right\Vert _{1}=\infty$
then implies zero Fisher information by the first part of the Theorem. On the
other hand, by Holder inequality with $p=1/\beta$ and $q=1/(1-\beta)$, for any
$0<\beta<1,$ and for all $b$ with $||b||_{\mathbf{H}}\leq1,$%
\begin{align*}
\left\vert \dot{\phi}(b)\right\vert  &  \leq\left(  \sum_{j=1}^{\infty}%
\lambda_{j}^{-2\beta}\left\langle r_{\phi},\varphi_{j}\right\rangle
_{_{\mathbf{H}}}^{2}\right)  ^{1/2}\left(  \sum_{j=1}^{\infty}\lambda
_{j}^{2\beta}\left\langle b,\varphi_{j}\right\rangle _{_{\mathbf{H}}}%
^{2}\right)  ^{1/2}\\
&  \leq\left\Vert r_{\phi}\right\Vert _{\beta}\left(  \sum_{j=1}^{\infty
}\lambda_{j}^{2\beta}\left\langle b,\varphi_{j}\right\rangle _{_{\mathbf{H}}%
}^{2\beta}\left\langle b,\varphi_{j}\right\rangle _{_{\mathbf{H}}}^{2-2\beta
}\right)  ^{1/2}\\
&  \leq\left\Vert r_{\phi}\right\Vert _{\beta}\left(  \sum_{j=1}^{\infty
}\lambda_{j}^{2}\left\langle b,\varphi_{j}\right\rangle _{_{\mathbf{H}}}%
^{2}\right)  ^{\beta/2}\left(  \sum_{j=1}^{\infty}\left\langle b,\varphi
_{j}\right\rangle _{_{\mathbf{H}}}^{2}\right)  ^{(1-\beta)/2}\\
&  \leq\left\Vert r_{\phi}\right\Vert _{\beta}||Sb||^{\beta}.
\end{align*}
Thus, for $\rho=1/\beta,$%
\[
\inf_{b\in\mathcal{B}_{\phi},||b||_{\mathbf{H}}\leq1}\frac{||Sb||^{2}}%
{[\dot{\phi}(b)]^{2\rho}}\geq\frac{1}{\left\Vert r_{\phi}\right\Vert _{\beta
}^{2\rho}}>0,
\]
and local identification follows. $\blacksquare$\bigskip

\noindent\textbf{Proof of Proposition \ref{Charactpar}}: For the functional
$\phi(\lambda)=\theta$ it holds $\mathcal{N}(\dot{\phi})=\{(b_{\theta}%
,b_{\eta}):b_{\theta}=0\}$. Then, by orthogonality
\begin{align*}
\mathcal{N}(S)  &  =\left\{  (b_{\theta},b_{\eta}):\int(\dot{l}_{\theta
}^{\prime}b_{\theta}+\dot{l}_{\eta}b_{\eta})^{2}d\mathbb{P}_{\theta_{0}%
,\eta_{0}}=0\right\} \\
&  =\left\{  (b_{\theta},b_{\eta}):\int(\tilde{l}_{\theta}^{\prime}b_{\theta
}+\Pi_{\overline{\mathcal{R}(\dot{l}_{\eta})}}\dot{l}_{\theta}^{\prime
}b_{\theta}+\dot{l}_{\eta}b_{\eta})^{2}d\mathbb{P}_{\theta_{0},\eta_{0}%
}=0\right\} \\
&  =\left\{  (b_{\theta},b_{\eta}):\int\left(  \tilde{l}_{\theta}^{\prime
}b_{\theta}\right)  ^{2}d\mathbb{P}_{\theta_{0},\eta_{0}}=0,\int%
(\Pi_{\overline{\mathcal{R}(\dot{l}_{\eta})}}\dot{l}_{\theta}^{\prime
}b_{\theta}+\dot{l}_{\eta}b_{\eta})^{2}d\mathbb{P}_{\theta_{0},\eta_{0}%
}=0\right\} \\
&  =\left\{  (b_{\theta},b_{\eta}):b_{\theta}^{\prime}\tilde{I}_{\theta
}b_{\theta}=0,\Pi_{\overline{\mathcal{R}(\dot{l}_{\eta})}}\dot{l}_{\theta
}^{\prime}b_{\theta}=-\dot{l}_{\eta}b_{\eta}\right\}  .
\end{align*}
Then, if $\tilde{I}_{\theta}>0,$ then we have $\mathcal{N}(S)\subset
\mathcal{N}(\dot{\phi})$. If $\tilde{I}_{\theta}=0,$ then there are two cases:
(1) $\dot{l}_{\theta}\in\mathcal{R}(\dot{l}_{\eta})$ and (2) $\dot{l}_{\theta
}\in\overline{\mathcal{R}(\dot{l}_{\eta})}\diagdown\mathcal{R}(\dot{l}_{\eta
}).$ In case (1) the identification condition $\mathcal{N}(S)\subset
\mathcal{N}(\dot{\phi})$ does not hold, as we can find $b_{\theta}\neq0$ such
that $\dot{l}_{\theta}b_{\theta}=-\dot{l}_{\eta}b_{\eta}$ (so $(b_{\theta
},b_{\eta})\in\mathcal{N}(S)$ but $(b_{\theta},b_{\eta})\notin\mathcal{N}%
(\dot{\phi})).$ In case (2) $\mathcal{N}(S)\subset\mathcal{N}(\dot{\phi})$
holds even though there is zero information for the parameter. Thus, if
$\dot{l}_{\theta}\in\overline{\mathcal{R}(\dot{l}_{\eta})}\diagdown
\mathcal{R}(\dot{l}_{\eta})$ we have \textquotedblleft irregular
identification\textquotedblright. $\blacksquare$\bigskip

\noindent\textbf{Proof of Theorem \ref{NonlinearPositive}}: Choose
$0<\varepsilon<I_{\theta,\rho}^{1/2},$ and $\delta>0$ such that for all
$\lambda\in\mathcal{B}_{\delta}(\lambda_{0})$ with $\theta\neq\theta_{0}$%
\begin{align}
\frac{\left\Vert \left(  f_{\lambda}-f_{\lambda_{0}}\right)  /f_{\lambda_{0}%
}-S(\lambda-\lambda_{0})\right\Vert }{\left\Vert S(\lambda-\lambda
_{0})\right\Vert }  &  =\frac{\left\Vert \left(  f_{\lambda}-f_{\lambda_{0}%
}\right)  /f_{\lambda_{0}}-S(\lambda-\lambda_{0})\right\Vert }{\left\vert
\theta-\theta_{0}\right\vert ^{\rho}}\frac{\left\vert \theta-\theta
_{0}\right\vert ^{\rho}}{\left\Vert S(\lambda-\lambda_{0})\right\Vert
}\nonumber\\
&  \leq\varepsilon\times I_{\theta,\rho}^{-1/2}\nonumber\\
&  <1, \label{ine1}%
\end{align}
where have used Assumption 2 and the definition of the generalized Fisher
information. The inequality (\ref{ine1}) implies that $\left\Vert \left(
f_{\lambda}-f_{\lambda_{0}}\right)  /f_{\lambda_{0}}\right\Vert \neq0,$ or
equivalently $f_{\lambda}\neq f_{\lambda_{0}}.$ That is, local identification
holds. $\blacksquare$\bigskip

\noindent\textbf{Proof of Proposition \ref{Ex1}}: Define $L:L_{1}(\pi)\mapsto
L_{1}(\mathbb{P})$ as
\[
Lb=\int f_{z/\alpha,\beta}(t_{1},t_{2})b(\alpha,\beta)d\pi(\alpha,\beta),
\]
where the conditional density of $Z$ given $(\alpha,\beta)$ is%
\[
f_{z/\alpha,\beta}(t_{1},t_{2})\propto f(t_{1};\alpha,\beta)f(t_{2}%
;\alpha,\beta),
\]
($\propto$ denotes equality up to multiplication by a normalizing constant)
and $f(t;\alpha,\beta)$ denotes the inverse Gaussian density%
\[
f(t;\alpha,\beta)\propto\frac{\beta}{t^{3/2}}e^{-\frac{\left(  \alpha
t-\beta\right)  ^{2}}{2t}}.
\]
Decompose $L$ as%
\begin{align*}
Lb  &  =\int_{\alpha\geq0}f_{z/\alpha,\beta}(t_{1},t_{2})b(\alpha,\beta
)d\pi(\alpha,\beta)+\int_{\alpha\leq0}f_{z/\alpha,\beta}(t_{1},t_{2}%
)b(\alpha,\beta)d\pi(\alpha,\beta)\\
&  \equiv L_{+}b+L_{-}b.
\end{align*}
Theorem 1 in Alvarez et al. (2016) and Proposition \ref{Identification} here
imply that both $L_{+}$ and $L_{-}$ are linear injective operators, and
therefore have inverses, $L_{+}^{-1}$ and $L_{-}^{-1},$ respectively. Define
the normalizing positive constant%
\[
C_{L}=\int f(t_{1};\alpha,\beta)f(t_{2};\alpha,\beta)\lambda_{0}(\alpha
,\beta)d\pi(\alpha,\beta)dt_{1}dt_{2}.
\]
Then, using that the inverse Gaussian satisfies
\[
f(t;\alpha,\beta)=e^{2\alpha\beta}f(t;-\alpha,\beta),
\]
it can be shown that%
\begin{align*}
L_{-}b  &  =C_{L}^{-1}\int_{\alpha\leq0}e^{4\alpha\beta}f(t_{1};-\alpha
,\beta)f(t_{2};-\alpha,\beta)b(\alpha,\beta)d\pi(\alpha,\beta)\\
&  =-C_{L}^{-1}\int_{\alpha\geq0}e^{-4\alpha\beta}f(t_{1};\alpha,\beta
)f(t_{2};\alpha,\beta)b(-\alpha,\beta)d\pi(\alpha,\beta)\\
&  =-L_{+}(e^{-4\alpha\beta}b(-\alpha,\beta)).
\end{align*}
Then, using these results, $b\in\mathcal{N}(L)$, i.e. $L_{+}b+L_{-}b=0$, is
equivalent to%
\begin{align*}
b(\alpha,\beta)  &  =-L_{+}^{-1}\left(  L_{-}b\right) \\
&  =L_{+}^{-1}L_{+}(e^{-4\alpha\beta}b(-\alpha,\beta))\\
&  =e^{-4\alpha\beta}b(-\alpha,\beta).
\end{align*}
This concludes the proof after noticing that $\mathcal{N}(S)=\mathcal{N}(L)$.
$\blacksquare$\bigskip

\noindent\textbf{Proof of Proposition \ref{PropABS}}: Let $f_{z,\alpha,\beta
}(t_{1},t_{2},\alpha,\beta)$ denote the joint density of $(Z,\alpha,\beta)$
wrt the product of the Lebesgue measure and $\pi,$ and $f_{\alpha,\beta
/z}(t_{1},t_{2})$ the corresponding conditional density of $(\alpha,\beta)$
given $Z.$ By the definition of $S$ with the reparametrization of Remark 1(ii)%
\begin{align*}
Sb  &  =\frac{1}{f_{\lambda_{0}}(t_{1},t_{2})}\int f_{z/\alpha,\beta}%
(t_{1},t_{2})b(\alpha,\beta)\lambda_{0}(\alpha,\beta)d\pi(\alpha,\beta)\\
&  =\frac{1}{f_{\lambda_{0}}(t_{1},t_{2})}\int f_{z,\alpha,\beta}(t_{1}%
,t_{2})b(\alpha,\beta)d\pi(\alpha,\beta)\\
&  =\int f_{\alpha,\beta/z}(t_{1},t_{2})b(\alpha,\beta)d\pi(\alpha,\beta)\\
&  =\mathbb{E}\left[  \left.  b(\alpha,\beta)\right\vert Z\right]  .
\end{align*}
By well-known properties of adjoint operators and substitution of
$f_{z/\alpha,\beta}(t_{1},t_{2}),$ we obtain
\begin{align*}
S^{\ast}g  &  =\mathbb{E}\left[  \left.  g(Z)\right\vert \alpha,\beta\right]
.\\
&  =\int_{\mathcal{T}^{2}}g(t_{1},t_{2})f_{z/\alpha,\beta}(t_{1},t_{2}%
)dt_{1}dt_{2}\\
&  =C\beta^{2}e^{2\alpha\beta}h(\alpha^{2},\beta^{2}),
\end{align*}
where%
\[
h(u,v)=\int_{\mathcal{T}^{2}}g(t_{1},t_{2})\frac{1}{t_{1}^{3/2}t_{2}^{3/2}%
}s(u,v;t_{1})s(u,v;t_{2})dt_{1}dt_{2}%
\]
and
\[
s(u,v;t)=\exp\left(  -\frac{ut}{2}-\frac{v}{2t}\right)  ,\text{ }%
t\in\mathcal{T},\text{ }(u,v)\in(0,\infty).
\]
We check that the conditions for an application of the Leibniz's rule hold.
These conditions are

\begin{description}
\item[1.] The partial derivative $\partial^{m}s(u,v;t_{1})s(u,v;t_{2}%
)/\partial^{m}u$ exists and is a continuous function on an open neighborhood
$B$ of $(u,v),$ for a.s. $(t_{1},t_{2})\in\mathcal{T}^{2}.$

\item[2.] There is a positive function $h_{m}(t_{1},t_{2})$ such that%
\begin{equation}
\sup_{(u,v)\in B}\left\vert \frac{\partial^{m}s(u,v;t_{1})s(u,v;t_{2}%
)}{\partial^{m}u}\right\vert \leq h_{m}(t_{1},t_{2}) \label{dom}%
\end{equation}
and%
\begin{equation}
\int_{\mathcal{T}^{2}}g(t_{1},t_{2})\frac{1}{t_{1}^{3/2}t_{2}^{3/2}}%
h_{m}(t_{1},t_{2})dt_{1}dt_{2}<\infty. \label{int}%
\end{equation}

\end{description}

\noindent Simple differentiation and induction show that for any integer
$m\geq0$%
\[
\frac{\partial^{m}s(u,v;t_{1})s(u,v;t_{2})}{\partial^{m}u}=2^{-m}%
(-1)^{m}(t_{1}+t_{2})^{m}s(u,v;t_{1})s(u,v;t_{2}).
\]
Therefore, by monotonicity we can find $u^{\ast}$ and $v^{\ast}$ such that
(\ref{dom}) holds with
\[
h_{m}(t_{1},t_{2})=2^{-m}(t_{1}+t_{2})^{m}s(u^{\ast},v^{\ast};t_{1})s(u^{\ast
},v^{\ast};t_{2}).
\]
Furthermore, by $\mathbb{E}\left[  \left.  g(Z)\right\vert \alpha
,\beta\right]  <\infty$ for all $\alpha$ and $\beta,$ and the boundedness of
$\mathcal{T}$, condition (\ref{int}) holds. The continuity of $h(u,v)$ is
proved using similar arguments (i.e. for $m=0)$. $\blacksquare$\bigskip

\noindent\textbf{Proof of Corollary \ref{PropABScor}}: By Theorem
\ref{MainTheoremS} part (i) it is sufficient to check that $r_{\phi}%
\notin\mathcal{R}(S^{\ast}),$ provided Assumption 1(i,ii,iv) holds. By
Proposition \ref{PropABS} the score operator is well-defined linear and
continuous. Furthermore, because heterogeneity is nonparametric and the
discussion before the Corollary, the corresponding functional $\dot{\phi}$ for
$\phi(\lambda_{0})=\mathbb{E}\left[  1(\alpha\leq\alpha_{0})1(\beta\leq
\beta_{0})\right]  $ is%
\[
\dot{\phi}(b)=\langle r_{\phi},b\rangle_{\mathbf{H}},
\]
where $\mathbf{H}=L_{2}^{0}(G_{0})$ and the Riesz representer $r_{\phi}%
\in\mathbf{H}$ is%
\[
r_{\phi}(\alpha,\beta)=1(\alpha\leq\alpha_{0})1(\beta\leq\beta_{0}%
)-G(\alpha_{0},\beta_{0}).
\]
Clearly, this function $r_{\phi}\notin\mathcal{R}(S^{\ast})$ and the result
follows. $\blacksquare$\bigskip

\noindent\textbf{Proof of Proposition \ref{WTP0}}: The score operator and the
functional satisfy Assumptions 1(i,ii,iv), since $r\in L_{2}(\pi)$. The
condition%
\[
r_{\phi}(w,x)=0\text{ }\pi-\text{a.s on }v_{\max}<w\leq w_{\max},
\]
then implies (\ref{1}). By Proposition \ref{Identification} and the discussion
of Assumption 1, see also Remark 2, identification follows. $\blacksquare
$\bigskip

\noindent\textbf{Proof of Proposition \ref{WTP1}}: To show that (\ref{genreg})
is sufficient for regular identification, note that (\ref{sol}) is a solution
to $r=S^{\ast}g$ with $g\in L_{2}.$ The see that $g\in L_{2}$, note%
\begin{align*}
\mathbb{E}\left[  g^{2}(Z)\right]   &  =\mathbb{E}\left[  g^{2}(1,V,X)G_{0}%
(V,X)+g^{2}(0,V,X)(1-G_{0}(V,X))\right] \\
&  \leq\mathbb{E}\left[  g^{2}(1,V,X)+g^{2}(0,V,X)\right] \\
&  \leq2\int_{0}^{v_{\max}}\int\left[  \frac{\partial r(v,x)}{\partial
v}\right]  ^{2}\frac{1}{f_{V/X=x}(v)}d\mu_{X}(x)dv<\infty.
\end{align*}
The necessity follows from the arguments prior to the statement of the
Proposition, i.e. from%
\[
\frac{\partial r(w,x)}{\partial w}=\left(  g(0,w,x)-g(1,w,x)\right)
f_{V/X=x}(w),
\]
and the fact that $g\in L_{2}.$ Conclude by an application of Theorem
\ref{MainTheoremS}(i). $\blacksquare$\bigskip

\noindent\textbf{Proof of Proposition \ref{WTP2}}: I first verify that median
WTP is identified. This is implied by continuity $G_{0}\left(  v\right)  $ and
the support condition $v_{\max}\geq\phi(\lambda_{0}),$ since by the
conditional independence%
\[
G_{0}\left(  w\right)  =\mathbb{E[E}\left[  \left.  Y\right\vert X,V=w\right]
]
\]
is identified on $\mathcal{S}_{V}$ and $G_{0}\left(  v_{\max}\right)  \geq1/2$
by monotonicity of cdfs. Standard arguments on quantiles, see Corollary 21.5
in van der Vaart (1998), yield that the median\ WTP has an influence function%
\[
r_{\phi}(w)=\frac{-\left\{  1(w<\phi(\lambda_{0}))-0.5\right\}  }%
{\mathbb{E}\left[  \lambda_{0}\left(  \phi(\lambda_{0}),X_{i}\right)  \right]
}.
\]
The discontinuity of the influence function then implies that identification
must be irregular, since $S^{\ast}g(w,x)$ is absolutely continuous in $w$.
$\blacksquare$\bigskip

\noindent\textbf{Proof of Proposition \ref{discount}}: By (\ref{A2dis})
Assumption 2 holds with $\rho=1$. By Theorem \ref{NonlinearPositive} it
remains to check $I_{\theta,1}>0.$ This is, however, equivalent to
$(1,0)\in\mathcal{R}(S^{\ast})$, or existence of $g\in L_{2}(\mu)$ such that%
\[
\langle A\eta_{0},g\rangle=1,\dot{l}_{\eta}^{\ast}g=0.
\]
If $\langle\eta_{0},g_{0}\rangle\neq0$ for $g_{0}$ one of the eigenfunctions
of $A^{\ast},$ we define $g=cg_{0}$ with $c=\left(  \theta_{0}\langle\eta
_{0},g_{0}\rangle\right)  ^{-1}.$ Note that, $\dot{l}_{\eta}^{\ast}g=c\dot
{l}_{\eta}^{\ast}g_{0}=0$ and $\langle A\eta_{0},g\rangle=c\theta_{0}%
\langle\eta_{0},g_{0}\rangle=1.$ $\blacksquare$\bigskip

\noindent\textbf{Proof of Proposition \ref{AARA}}: The functional $\chi
(\dot{u}_{0}),$ although nonlinear, is concave and differentiable, with
pathwise derivative
\[
\dot{\chi}(b)=\mathbb{E}\left[  b(C_{t}^{\ast})\frac{d(C_{t}^{\ast})}{\dot
{u}_{0}(C_{t}^{\ast})}\right]  ,
\]
which implies, with $\dot{l}_{\eta}\eta=\theta_{0}A\eta(c)-\eta(c),$
\begin{align*}
\varpi(\epsilon)  &  =\sup_{\left\Vert \dot{u}_{1}-\dot{u}_{0}\right\Vert
_{2,\mu^{\ast}}\leq\delta,||\dot{l}_{\eta}L(\dot{u}_{1}-\dot{u}_{0})||_{2,\mu
}\leq\epsilon}\left\vert \chi(\dot{u}_{1})-\chi(\dot{u}_{0})\right\vert \\
&  \leq\delta\sup_{\left\Vert \dot{u}_{1}-\dot{u}_{0}\right\Vert _{2,\mu
^{\ast}}\leq\delta,||\dot{l}_{\eta}L(\dot{u}_{1}-\dot{u}_{0})||_{2,\mu}%
\leq\epsilon}\left\vert \dot{\chi}(\dot{u}_{1}-\dot{u}_{0})\right\vert \\
&  =\delta\sup_{\left\Vert \dot{u}_{1}-\dot{u}_{0}\right\Vert _{2,\mu^{\ast}%
}\leq\delta,||\dot{l}_{\eta}L(\dot{u}_{1}-\dot{u}_{0})||_{2,\mu}\leq\epsilon
}\left\vert \mathbb{E}\left[  \dot{l}_{\eta}L(\dot{u}_{1}(C_{t})-\dot{u}%
_{0}(C_{t}))g(C_{t})\right]  \right\vert \\
&  \leq\delta||g||_{2,\mu}\epsilon,
\end{align*}
where the first inequality uses concavity, the last equality uses
(\ref{AARAreg}), so
\begin{align*}
\dot{\chi}(\dot{u}_{1}-\dot{u}_{0})  &  =\mathbb{E}\left[  (\dot{u}_{1}%
-\dot{u}_{0})(C_{t}^{\ast})r_{\chi}(C_{t}^{\ast})\right] \\
&  =\mathbb{E}\left[  (\dot{u}_{1}-\dot{u}_{0})(C_{t}^{\ast})\dot{l}_{\eta
}^{\ast}g(C_{t})\right] \\
&  =\mathbb{E}\left[  \dot{l}_{\eta}L(\dot{u}_{1}(C_{t})-\dot{u}_{0}%
(C_{t}))g(C_{t})\right]  ,
\end{align*}
and the last inequality follows by Cauchy-Schwarz. Lemma \ref{LemmaMC} in the
Supplemental Appendix then implies local identification of $\chi(\dot{u}_{0}%
)$. It remains to show the regularity, but this follows from Cauchy-Schwarz,
since%
\[
\inf_{\dot{u}:\dot{\chi}(\dot{u})\neq0}\frac{||\dot{l}_{\eta}L\dot{u}||^{2}%
}{\left\vert \dot{\chi}(\dot{u})\right\vert ^{2}}\geq\frac{\left\Vert \dot
{l}_{\eta}L\dot{u}\right\Vert ^{2}}{||g||_{2,\mu}^{2}\left\Vert \dot{l}_{\eta
}L\dot{u}\right\Vert ^{2}}=\frac{1}{||g||_{2,\mu}^{2}}>0.
\]
$\blacksquare$

\noindent\textbf{Proof of Lemma \ref{LemmaP1}}: By the change of variables
with $c=\exp(z)$ and $c^{\ast}=\exp(\tau),$ and multiplying both sides by
$\exp(-\tau)$ the integral equation
\[
r_{\chi}(c^{\ast})=\int f_{\varepsilon}(c/c^{\ast})w(c)dc
\]
is transformed into a convolution-type problem
\begin{equation}
x(\tau)=\int K(z-\tau)y(z)dz, \label{conv}%
\end{equation}
where $x(\tau)=\exp(-\tau)r_{\chi}(\exp(\tau))$, $K(u)=\exp(u)f_{\varepsilon
}(\exp(u))$ and $y(z)=w(\exp(z)).$ By Polyanin and Manzhirov (2008, p. 285) if
$x\in L_{2}(\mathbb{R})$, a necessary and sufficient condition for existence
of $y\in L_{2}(\mathbb{R})$ satisfying (\ref{conv}) is $\hat{x}/\hat{K}\in
L_{2}(\mathbb{R})$. The solution is given by%
\[
y(z)=\frac{1}{2\pi}\operatorname{Re}\int e^{itz}\frac{\hat{x}(t)}{\hat{K}%
(t)}dt,
\]
and in terms of $w$,%
\[
w(c)=\frac{1}{2\pi}\operatorname{Re}\int e^{it\log(c)}\frac{\hat{x}(t)}%
{\hat{K}(t)}dt.
\]
Note this solution is also in $L_{2}(\mu)$ if $\hat{x}/\hat{K}\in
L_{2}(\mathbb{R})$, since, by a change of variables and Fubini,
\begin{align*}
\mathbb{E}\left[  \left\vert w(C_{t})\right\vert ^{2}\right]   &  =\frac
{1}{\left(  2\pi\right)  ^{2}}\left\vert \int\int\frac{\hat{x}(t)}{\hat{K}%
(t)}\frac{\hat{x}(s)}{\hat{K}(s)}\mathbb{E[}e^{i(t-s)\log C_{t}}%
]dtds\right\vert \\
&  \leq\frac{1}{\left(  2\pi\right)  ^{2}}\int\left\vert \frac{\hat{x}%
(t)}{\hat{K}(t)}\right\vert ^{2}dt<\infty.
\end{align*}
\newpage

\section{Supplemental Appendix}

\label{Supplement}

I will extensively use basic results from operator theory and Hilbert spaces
in this Supplemental Material. See Carrasco, Florens and Renault (2007) for an
excellent review of these results. This Appendix is organized as follows.
Section \ref{Linearwithout} establishes sufficient conditions for local
irregular identification in models linear in nuisance parameters. Section
\ref{Nuisance} characterizes identification of linear continuous functionals
of nuisance parameters in semiparametric models. Section \ref{Nonlinear}
establishes sufficient conditions for identification in general nonlinear models.

\subsection{Models Linear in Nuisance Parameters}

\label{Linearwithout}

Define the nuisance score operator
\begin{equation}
\dot{l}_{\eta(\theta)}b_{\eta}=\frac{f_{\theta,\eta_{0}+b_{\eta}}%
-f_{\theta,\eta_{0}}}{f_{\theta_{0},\eta_{0}}}, \label{aso}%
\end{equation}
and the (negative) approximated score for $\theta$ as%
\[
s_{\theta}=\frac{f_{\theta_{0},\eta_{0}}-f_{\theta,\eta_{0}}}{f_{\theta
_{0},\eta_{0}}}.
\]
I drop the dependence on $\theta_{0}$ and denote $\dot{l}_{\eta}\equiv\dot
{l}_{\eta(\theta_{0})}.$ Define the (negative) approximated efficient score
$\tilde{s}_{\theta}:=s_{\theta}-\Pi_{\overline{\mathcal{R}(\dot{l}%
_{\eta(\theta)})}}s_{\theta},$ and the approximated Fisher Information%
\[
G(\theta)=||\tilde{s}_{\theta}||^{2}.
\]
Let $\Psi$ be the class of measurable functions $\psi:[0,\infty
)\longrightarrow\lbrack0,\infty)$ that are increasing, right continuous at $0$
and with $\psi(0)=0.$ Then, consider the following assumption. \bigskip

\noindent\textbf{Assumption D}: (i) The map $\dot{l}_{\eta(\theta)}:T(\eta
_{0})\subseteq\mathcal{H}\mapsto L_{2}$ is linear for each $\theta$ in a
neighborhood of $\theta_{0}$ (ii) there exists a positive constant $C$ such
that $G(\theta)>C\psi(\left\vert \theta-\theta_{0}\right\vert ^{2})$ in a
neighborhood of $\theta_{0},$ where $\psi\in\Psi.$\bigskip

\noindent Assumption D(i) holds for many models of interest. Assumption D(ii)
follows from conditions on the derivative of $G(\theta)$ at $\theta_{0}.$ For
example, if $G(\theta)$ is differentiable at $\theta_{0}$ with full rank
derivative at $\theta_{0}$, then Assumption D(ii) holds with $\psi
(\epsilon)=\epsilon.$ This corresponds to the case of regular local
identification. A necessary condition for Assumption D(ii) is that
$\mathcal{N}(\dot{l}_{\eta(\theta)}^{\ast})\neq0,$ since otherwise
$G(\theta)=0$.\bigskip

\begin{theorem}
\label{Semiparametric}Let Assumption D hold. Then, $\theta$ is locally
identified at $\theta_{0}$.
\end{theorem}

\noindent\textbf{Proof of Theorem \ref{Semiparametric}}: Write%
\begin{align*}
\frac{f_{\theta,\eta}-f_{\theta_{0},\eta_{0}}}{f_{\theta_{0},\eta_{0}}}  &
=\frac{f_{\theta,\eta}-f_{\theta,\eta_{0}}}{f_{\theta_{0},\eta_{0}}}%
-\frac{f_{\theta_{0},\eta_{0}}-f_{\theta,\eta_{0}}}{f_{\theta_{0},\eta_{0}}}\\
&  =\dot{l}_{\eta(\theta)}b_{\eta}-s_{\theta}.
\end{align*}
Note that by standard least squares theory for all $b_{\eta}\in T(\eta_{0}),$
and all $\theta$ in a neighborhood of $\theta_{0},$
\begin{align*}
||\dot{l}_{\eta(\theta)}b_{\eta}-s_{\theta}||^{2}  &  \geq||\Pi_{\overline
{\mathcal{R}(\dot{l}_{\eta(\theta)})}}s_{\theta}-s_{\theta}||^{2}\\
&  >C\psi(\left\vert \theta-\theta_{0}\right\vert ^{2}).
\end{align*}
This inequality implies local identification. $\blacksquare$\bigskip

\subsection{Functionals of Nuisance Parameters in Semiparametric Models}

\label{Nuisance}

\noindent Let $\chi:\mathcal{H}\mapsto\mathbb{R}$ be a linear continuous
functional, and let $r_{\chi}\in T(\eta_{0})\subset\mathcal{H}$ be such that
for all $b_{\eta}\in T(\eta_{0}),$
\[
\chi(b_{\eta})=\langle b_{\eta},r_{\chi}\rangle_{\mathcal{H}}.
\]
To give a general result, I allow for $\theta$ to be infinite-dimensional, and
ask the question: When lack of identification of one parameter, here $\theta,$
does not have an effect, at least locally, on identification on another
parameter $\chi(\eta)$?

A similar characterization to that of Proposition \ref{Charactpar} is obtained
for $\phi(\lambda)=\chi(\eta)$, allowing for singular information for both
$\theta$ and the functional $\phi(\lambda)=\chi(\eta)$. Define the operator%
\[
A_{\eta\theta}=\left(  \dot{l}_{\eta}^{\ast}\dot{l}_{\eta}\right)  ^{-}\dot
{l}_{\eta}^{\ast}\dot{l}_{\theta},
\]
where $B^{-}$ denotes the generalized Moore-Penrose inverse of $B.$\bigskip

\begin{proposition}
\label{Charactnui}For the functional $\phi(\lambda)=\chi(\eta)\in\mathbb{R}$:
(i) if $\mathcal{R}(\dot{l}_{\theta})\cap\mathcal{R}(\dot{l}_{\eta})=\{0\},$
then $\mathcal{N}(S)\subset\mathcal{N}(\dot{\phi})$ holds iff $r_{\chi}%
\in\overline{\mathcal{R}(\dot{l}_{\eta}^{\ast})};$ (ii) if $\mathcal{R}%
(\dot{l}_{\theta})\cap\mathcal{R}(\dot{l}_{\eta})\neq\{0\},$ then
$\mathcal{N}(S)\subset\mathcal{N}(\dot{\phi})$ holds if $r_{\chi}\in
\overline{\mathcal{R}(\dot{l}_{\eta}^{\ast})}\cap\mathcal{N}(A_{\eta\theta
}^{\ast})$.\ 
\end{proposition}

\noindent\textbf{Proof of Proposition \ref{Charactnui}}: Note that for the
functional $\phi(\lambda)=\chi(\eta),$ where $\chi:H\mapsto\mathbb{R}$ is a
linear continuous functional with
\[
\chi(b_{\eta})=\langle b_{\eta},r_{\chi}\rangle_{H},
\]
it holds that $\mathcal{N}(\dot{\phi})=\{(b_{\theta},b_{\eta}):\langle
b_{\eta},r_{\chi}\rangle_{H}=0\}$. Therefore, by the proof of Proposition
\ref{Charactpar} (which is also valid for infinite-dimensional $\theta,$ with
$\tilde{I}_{\theta}$ interpreted as an operator), $\mathcal{N}(S)\subset
\mathcal{N}(\dot{\phi})$ iff $b_{\theta}^{\prime}\tilde{I}_{\theta}b_{\theta
}=0$ and $\Pi_{\overline{\mathcal{R}(\dot{l}_{\eta})}}\dot{l}_{\theta}%
^{\prime}b_{\theta}=-\dot{l}_{\eta}b_{\eta}$ implies $\langle b_{\eta}%
,r_{\chi}\rangle_{H}=0.$ If $\tilde{I}_{\theta}$ is positive definite, then
$(b_{\theta},b_{\eta})\in\mathcal{N}(S)$ iff $b_{\theta}=0$ and $0=\dot
{l}_{\eta}b_{\eta}.$ Therefore, $(b_{\theta},b_{\eta})\in\mathcal{N}(\dot
{\phi})$ iff $\mathcal{N}(\dot{l}_{\eta})\subset\mathcal{N}(\chi),$ which is
equivalent to $r_{\chi}\in\overline{\mathcal{R}(\dot{l}_{\eta}^{\ast}).}$ If
$\tilde{I}_{\theta}$ is semi-positive definite, there are two cases (i)
$\mathcal{R}(\dot{l}_{\theta})\cap\mathcal{R}(\dot{l}_{\eta})\neq\{0\}$ and
(ii) $\mathcal{R}(\dot{l}_{\theta})\subset\overline{\mathcal{R}(\dot{l}_{\eta
})}\diagdown\mathcal{R}(\dot{l}_{\eta}).$ In case (i), $\dot{l}_{\theta
}b_{\theta}=-\dot{l}_{\eta}b_{\eta},$ and for all such $b_{\eta}$ it must hold
that $\langle b_{\eta},r_{\chi}\rangle_{H}=0.$ All the solutions of $\dot
{l}_{\theta}b_{\theta}=-\dot{l}_{\eta}b_{\eta}$ can be written as $b_{\eta
}=\mathcal{N}(\dot{l}_{\eta})-A_{\eta\theta}b_{\theta}.$ Thus, the
orthogonality $\langle b_{\eta},r_{\chi}\rangle_{H}=0$ holds if $r_{\chi}%
\in\overline{\mathcal{R}(\dot{l}_{\eta}^{\ast})}\cap\mathcal{N}(A_{\eta\theta
}^{\ast}).$ In case (ii) $0=\dot{l}_{\eta}b_{\eta}$ must imply that
$(b_{\theta},b_{\eta})\in\mathcal{N}(\dot{\phi})$, which holds if
$\mathcal{N}(\dot{l}_{\eta})\subset\mathcal{N}(\chi)\ $or equivalently
$r_{\chi}\in\overline{\mathcal{R}(\dot{l}_{\eta}^{\ast}).}$ Therefore, if
$\mathcal{R}(\dot{l}_{\theta})\cap\mathcal{R}(\dot{l}_{\eta})=\{0\}$
($\tilde{I}_{\theta}\ $is positive definite or case (ii) above) then
$\mathcal{N}(S)\subset\mathcal{N}(\dot{\phi})$ holds iff $r_{\chi}\in
\overline{\mathcal{R}(\dot{l}_{\eta}^{\ast})};$ (ii) if $\mathcal{R}(\dot
{l}_{\theta})\cap\mathcal{R}(\dot{l}_{\eta})\neq\{0\}$ (case (i) above) then
$\mathcal{N}(S)\subset\mathcal{N}(\dot{\phi})$ holds if $\overline
{\mathcal{R}(\dot{l}_{\eta}^{\ast})}\cap\mathcal{N}(A_{\eta\theta}^{\ast}%
)$.\ $\blacksquare$\bigskip

\begin{remark}
The conditions for local identification of $\chi(\eta_{0})$ depend on whether
$\theta_{0}$ is locally identified or not. The case (ii) corresponds to the
situation of local unidentification of $\theta_{0},$ and it is shown that
despite this lack of local identification of $\theta_{0},$ $\chi(\eta_{0})$
might still be locally identified. To interpret the result, one can think of
$r_{\chi}\in\overline{\mathcal{R}(\dot{l}_{\eta}^{\ast})}$ as the
identification condition for $\chi(\eta_{0})$ that would be needed if
$\theta_{0}$ was known. If $\theta_{0}$ is not known, but is identified, one
can treat it as known for the purpose of identifying $\chi(\eta_{0}).$
However, if $\theta_{0}$ is not identified, an additional condition must be
met to avoid the lack of identification of $\theta_{0}$ to spread out to
$\chi(\eta_{0}).$ Technically, this condition is that for all $b=(b_{\theta
},b_{\eta})$ such that $\dot{l}_{\theta}b_{\theta}=-\dot{l}_{\eta}b_{\eta}$
(these $b^{\prime}s$ are directions that lead to zero nonparametric
information), it must hold that $\langle b_{\eta},r_{\chi}\rangle_{H}=0.$
Under $r_{\chi}\in\overline{\mathcal{R}(\dot{l}_{\eta}^{\ast})},$ a simple
condition for this orthogonality is $r_{\chi}\in\mathcal{N}(A_{\eta\theta
}^{\ast})$.
\end{remark}

\begin{remark}
In both cases $r_{\chi}\in\overline{\mathcal{R}(\dot{l}_{\eta}^{\ast}%
)}\diagdown\mathcal{R}(\dot{l}_{\eta}^{\ast})$ corresponds to the case of zero
information for $\phi(\lambda)=\chi(\eta)$ at $\phi(\lambda_{0})=\chi(\eta
_{0}).$ Regular identification of $\chi(\eta)$ in case (ii) requires that for
all $r_{\chi}^{\ast}$ that solve $r_{\chi}=\dot{l}_{\eta}^{\ast}r_{\chi}%
^{\ast}$ it holds that $r_{\chi}^{\ast}\in\mathcal{N}(\dot{l}_{\theta}^{\ast
})$. Under this condition, lack of identification of $\theta_{0}$ does not
affect regular identification of $\chi(\eta_{0}).$
\end{remark}

\noindent Van der Vaart (1991) has shown that a positive information of
$\chi(\eta_{0})$ is equivalent to $r_{\chi}\in\mathcal{R}(\dot{l}_{\eta}%
^{\ast})$ when $\theta_{0}$ is locally regularly identified and $\eta_{0}$ is
identified. Proposition \ref{Charactnui} characterizes local regular and
irregular identification of $\chi(\eta_{0}),$ allowing for $\theta_{0}$ to be
locally regular or irregularly identified, or even unidentified. The results
of Proposition \ref{Charactnui} are applied to measures of risk aversion in
Example 3 on the Euler Equation.

\subsection{General Nonlinear Models}

\label{Nonlinear}

The following modulus of continuity is shown to be useful for the study of
identification%
\begin{equation}
\varpi(\epsilon)=\sup_{\lambda\in\mathcal{B}_{\delta}(\lambda_{0}):||\left(
f_{\lambda}-f_{\lambda_{0}}\right)  f_{\lambda_{0}}^{-1}||\leq\epsilon
}\left\vert \phi(\lambda)-\phi(\lambda_{0})\right\vert . \label{modulus}%
\end{equation}
I drop the dependence of $\varpi(\epsilon)$ on $\delta$ for simplicity of
notation. Lemma \ref{LemmaMC} below shows that $\varpi(\epsilon)\downarrow0$
as $\epsilon\downarrow0$ is sufficient for local identification of
$\phi(\lambda_{0})$. A related modulus of continuity was introduced in Donoho
and Liu (1987) for the purpose of obtaining bounds on the optimal rate of
convergence for functionals of a density (they assume identification and use
the Hellinger metric). Using $||\left(  f_{\lambda}-f_{\lambda_{0}}\right)
f_{\lambda_{0}}^{-1}||$ is convenient because we can exploit simultaneously
the linearity of certain models and the Hilbert space structure. \bigskip

\noindent\textbf{Lemma }\label{LemmaMC} \textit{If there exists }$\delta
>0$\textit{ such that }$\varpi(\epsilon)\rightarrow0$\textit{ as }%
$\epsilon\rightarrow0,$\textit{ then }$\phi(\lambda_{0})$\textit{ is locally
identified.}

\noindent\textbf{Proof of Lemma \ref{LemmaMC}}: Suppose that $\phi(\lambda
_{0})$ is not locally identified. Then, for all $\delta>0$, we can find a
$\lambda^{\ast}\in\Lambda_{\delta}(\lambda_{0})$ such that $\left\Vert \left(
f_{\lambda^{\ast}}-f_{\lambda_{0}}\right)  /f_{\lambda_{0}}\right\Vert =0$ and
$\phi(\lambda^{\ast})\neq\phi(\lambda_{0}),$ and therefore, for all
$\epsilon>0,$
\[
\varpi(\epsilon)\geq\left\vert \phi(\lambda^{\ast})-\phi(\lambda
_{0})\right\vert >0,
\]
showing that $\varpi(\epsilon)$ does not converge to zero as $\epsilon
\rightarrow0.$ $\blacksquare$\bigskip

\noindent The following result provides a general local identification result.
Recall $\Psi$ is the class of measurable functions $\psi:[0,\infty
)\longrightarrow\lbrack0,\infty)$ that are increasing, right continuous at $0$
and with $\psi(0)=0.$ \bigskip

\noindent\textbf{Assumption N}:\textit{ For all }$\varepsilon>0,$\textit{
there exists }$\delta>0,$\textit{ }$\psi_{1},\psi_{2}\in\Psi$\textit{, and a
continuous linear operator }$S:T(\lambda_{0})\subseteq H\mapsto L_{2}%
,$\textit{ such that for all }$\lambda=(\theta,\eta)\in\mathcal{B}_{\delta
}(\lambda_{0}),$

\begin{description}
\item[(i)]
\[
\left\Vert \left(  f_{\lambda}-f_{\lambda_{0}}\right)  /f_{\lambda_{0}%
}-S(\lambda-\lambda_{0})\right\Vert <\varepsilon\psi_{1}\left(  \left\Vert
\lambda-\lambda_{0}\right\Vert _{\mathbf{H}}\right)  ;
\]

\item[(ii)]
\[
\left\vert \phi(\lambda)-\phi(\lambda_{0})\right\vert \leq\psi_{2}\left(
\left\Vert \lambda-\lambda_{0}\right\Vert _{\mathbf{H}}\right)  ;\text{
\textit{and}}%
\]

\item[(iii)]
\[
\inf_{\lambda\in\mathcal{B}_{\delta}(\lambda_{0})}\frac{||S(\lambda
-\lambda_{0})||}{\psi_{1}\left(  \left\Vert \lambda-\lambda_{0}\right\Vert
_{\mathbf{H}}\right)  }>0.
\]
\bigskip
\end{description}

\noindent Assumption N(i) and N(ii) are mild smoothness conditions that often
hold in applications. Condition N(iii) is a positive nonparametric generalized
information condition. Then, I have the following

\begin{theorem}
\label{GeneralNonlinear}Let Assumption N hold. Then, $\phi(\lambda)$ is
locally identified at $\phi(\lambda_{0})$.
\end{theorem}

\noindent\textbf{Proof of Theorem \ref{GeneralNonlinear}}: Assumptions N(i-ii)
imply that if $||\left(  f_{\lambda}-f_{\lambda_{0}}\right)  f_{\lambda_{0}%
}^{-1}||\leq\epsilon$ then we can find a positive constant $C$ and
$0<\varepsilon<C$ such that for all $\lambda=(\theta,\eta)\in\mathcal{B}%
_{\delta}(\lambda_{0}),$
\[
C\psi_{1}\left(  \left\Vert \lambda-\lambda_{0}\right\Vert _{\mathbf{H}%
}\right)  \leq\left\Vert S(\lambda-\lambda_{0})\right\Vert \leq\varepsilon
\psi_{1}\left(  \left\Vert \lambda-\lambda_{0}\right\Vert _{\mathbf{H}%
}\right)  +\epsilon,
\]
which in turn implies%
\[
\psi_{1}\left(  \left\Vert \lambda-\lambda_{0}\right\Vert _{\mathbf{H}%
}\right)  \leq\frac{\epsilon}{C-\varepsilon}.
\]
Hence, by Assumption N(ii)%
\begin{align*}
\varpi(\epsilon) &  =\sup_{\lambda\in\mathcal{B}_{\delta}(\lambda
_{0}):||\left(  f_{\lambda}-f_{\lambda_{0}}\right)  f_{\lambda_{0}}^{-1}%
||\leq\epsilon}\left\vert \phi(\lambda)-\phi(\lambda_{0})\right\vert ,\\
&  \leq\sup_{\lambda\in\mathcal{B}_{\delta}(\lambda_{0}):\psi_{1}\left(
\left\Vert \lambda-\lambda_{0}\right\Vert _{\mathbf{H}}\right)  \leq
\frac{\epsilon}{C-\varepsilon}}\psi_{2}\left(  \left\Vert \lambda-\lambda
_{0}\right\Vert _{\mathbf{H}}\right)  \\
&  \leq\psi_{2}\left(  \psi_{1}^{-1}\left(  \frac{\epsilon}{C-\varepsilon
}\right)  \right)  \\
&  \rightarrow0\text{ as }\epsilon\rightarrow0.
\end{align*}
Thus, the Theorem follows from Lemma \ref{LemmaMC}. $\blacksquare$

\subsubsection{A Counterexample}

\label{counterexample}

I provide a counterexample, building on that given in Chen et al. (2014, pg.
791), that shows that regular identification is not equivalent to $I_{\phi}>0$
in general (and hence to Van der Vaart's (1991) diffentiability condition).
Let $\lambda=(\lambda_{1},\lambda_{2},...)$ be a sequence of real numbers. Let
$(p_{1},p_{2},...)$ be probabilities, $p_{j}>0,$ $\sum_{j=1}^{\infty}p_{j}=1.$
Let $f(x)$ be a twice continuously differentiable function of a scalar $x$
that is bounded with bounded second derivative. Suppose $f(x)=0$ if and only
if $x\in\{0,1\}$ and $\partial f(0)/\partial x=1$. Let $m(\lambda
)=(f(\lambda_{1}),f(\lambda_{2}),...)$ also be a sequence with $\left\Vert
m(\lambda)\right\Vert ^{2}=\sum_{j=1}^{\infty}p_{j}f^{2}(\lambda_{j})$. Then,
for $\left\Vert \lambda\right\Vert _{\Lambda}=\left(  \sum_{j=1}^{\infty}%
p_{j}\lambda_{j}^{4}\right)  ^{1/4}$ the mapping $m$ is Frechet differentiable
at $\lambda_{0}=0$ with derivative $Sb=b$, but $\lambda_{0}=0$ is not locally
identified (Chen et al. 2014).

Consider the nonlinear functional%
\[
\phi(\lambda)=\sum_{j=1}^{\infty}f(\lambda_{j})p_{j}.
\]
This functional has a derivative at $\lambda_{0}=0$ given by%
\[
\dot{\phi}(b)=\sum_{j=1}^{\infty}b_{j}p_{j},
\]
and by Cauchy-Schwarz
\begin{align*}
\left\vert \dot{\phi}(b)\right\vert ^{2}  &  \leq\left(  \sum_{j=1}^{\infty
}b_{j}^{2}p_{j}\right) \\
&  =\left\Vert Sb\right\Vert ^{2}.
\end{align*}
Hence, $I_{\phi}\geq1>0.$ However, the functional is not identified, since
$\phi(\alpha^{k})=0=\phi(0),$ where $\alpha^{k}=(0,..,0,1,1,1...)$ has zeros
in the first $k$ positions and a one everywhere else. \newpage%

%TCIMACRO{\TeXButton{TeX field}{\begin{thebibliography}{99}}}%
%BeginExpansion
%
%EndExpansion

\end{document}